\documentclass{pspum-l}
\usepackage{enumerate}
\usepackage{graphicx}
\usepackage[ansinew]{inputenc}
\usepackage{natbib}  

\DeclareRobustCommand{\SkipTocEntry}[4]{}

\def\platz{\par\medbreak}

   {\end{trivlist}\platz}

\newtheorem*{main problem}{Main Problem}
\newtheorem*{Qgraph problem}{Quantum Graph Problem}

\newcommand\calD{{\mathcal D}}

\newcommand\calN{{\mathcal N}}

\newcommand\calX{{\mathcal X}}

\newcommand\diag{\operatorname{diag}}

\newcommand{\Dom}{\operatorname{Dom}}

\newcommand\Ran{\operatorname{Ran}}

\newcommand{\Span}{\operatorname{span}} 

\newcommand\vol{\operatorname{vol}}

\newcommand {\R}{\mathbb{R}}

\newcommand {\Z}{\mathbb{Z}}
\newcommand {\C}{\mathbb{C}}
\newcommand {\N}{\mathbb{N}}

\newcommand {\grad}{\nabla}

\newcommand{\Xinfty}{X^\infty}

\newcommand{\ubar}{\overline{u}}

\newcommand{\ave}{\operatorname{ave}}
\newcommand{\locmax}{\operatorname{loc.max.}}
\newcommand{\xbar}{\overline{x}}
\newcommand{\ybar}{\overline{y}}
\newcommand{\lambdabar}{\overline{\lambda}}
\newcommand{\Omegatilde}{\tilde{\Omega}}
\newcommand{\PDO}{$\Psi$DO}
\newcommand{\calXbar}{\overline{\calX}}

\theoremstyle{definition}

\theoremstyle{remark}

\numberwithin{equation}{section}

\begin{document}
\title[Thin tubes survey]{Thin tubes in mathematical physics, global analysis and spectral geometry}

\author{Daniel Grieser}
\address{Institut f\"ur Mathematik, Carl von Ossietzky Universit\"at Oldenburg,
D-26111 Oldenburg}
\email{grieser@mathematik.uni-oldenburg.de}
\keywords{Quantum graph, cylindrical end, gluing formula, nodal line, scattering}
\subjclass[2000]{Primary 58-02 
                         35-02 
                         81-02 
                         }
\date{February 11, 2008}
\begin{abstract}
A thin tube is an $n$-dimensional space which is very thin in $n-1$ directions, compared to the remaining direction, for example the $\varepsilon$-neighborhood of a curve or an embedded graph in $\R^n$ for small $\varepsilon$. The Laplacian on thin tubes and related operators have been studied in various contexts, with different goals but overlapping techniques. In this survey we explain some of these contexts, methods and results, hoping to encourage more interaction between the disciplines mentioned in the title.
\end{abstract}
\maketitle

\tableofcontents

\section{Problem}\label{sec problem}
A thin tube is, roughly speaking, a metric space which is close to a one-dimensional space, for example a curve. Often, but not always, it is useful to think of a thin tube as a member of a family $(\Omega_\varepsilon)_{0<\varepsilon<\varepsilon_0}$ where $\Omega_\varepsilon$ is built from a one-dimensional space $\Omega_0$  by adding transverse dimensions scaled to size $\varepsilon$ in a prescribed way. Thus, $\Omega_0$ can be thought of as limit (as $\varepsilon\to 0$) or model for the family $(\Omega_\varepsilon)_\varepsilon$.

\begin{figure}
\includegraphics[scale=0.8]{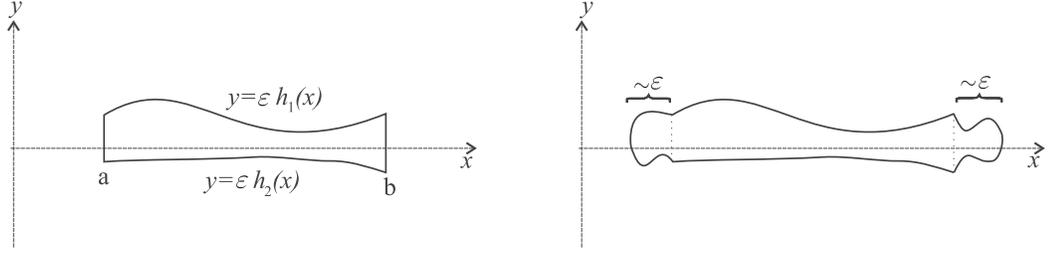}
\caption{Tube without ends and tube with ends}
\label{figure tubes}
\end{figure}

We will make this precise in various settings. We will consider two kinds of tubes: Tubes without ends and tubes with ends. To explain this, we define these first in the simple setting of two-dimensional tubes modeled on an interval. See Figure \ref{figure tubes}.
\medskip

\noindent{\bf Tubes without ends:}
Let $I\subset\R$ be a closed interval and let $h_1,h_2:I\to\R$ be continuous functions satisfying $h_1(x)<h_2(x)$ for all $x\in I$. For $\varepsilon>0$ the {tube without ends} is
\begin{equation}
\label{eqn def Omega1}
\Omega^{\text{without ends}}_\varepsilon = \{(x,y)\in\R^2:\ x\in I,\ \varepsilon h_1(x)\leq y\leq \varepsilon h_2(x)\}.
\end{equation}

\noindent{\bf Tubes with ends:}
In addition to the previous data, let $X_l$, $X_r$ be two bounded closed sets in $\R^2$ with non-empty interior such that $X_l$ is contained in the left half plane $\{x\leq 0\}$ and $X_r$ in the right half plane $\{x\geq 0\}$. A tube with ends is obtained by attaching $X_l$ and $X_r$, scaled by a factor $\varepsilon$, to the left and right of a tube without ends. Thus, for $I=[a,b]$
\begin{equation}
\label{eqn def Omega2}
\Omega^{\text{with ends}}_\varepsilon = \Omega^{\text{without ends}}_\varepsilon \cup [(a,0)+\varepsilon X_l] \cup [(b,0)+\varepsilon X_r].
\end{equation}
where $a+\varepsilon X:=\{(a+\varepsilon x,\varepsilon y):\ (x,y)\in X\}$.  We assume $X_l,X_r$ are such that $\Omega^{\text{with ends}}_\varepsilon$ is connected. Note that the ends are scaled by $\varepsilon$ in all directions while the tube is only scaled in the $y$-direction. This will make a substantial difference in the analysis.

We simply write $\Omega_\varepsilon$ if we do not want to specify which kind of tube we consider.
\medskip

The class of problems we will discuss concern the Laplace operator, $\Delta=-\left(\frac{\partial^2}{\partial x^2} + \frac{\partial^2}{\partial y^2}\right)$ on $\Omega_\varepsilon$. We will be mostly interested in the spectrum, that is, the solutions of
\begin{equation}
\label{eqn eigenvalues}
\Delta u = \lambda u,\quad u:\Omega_\varepsilon \to \C,\ \lambda\in\R
\end{equation}
which in addition satisfy boundary conditions, for example Dirichlet ($u_{|\partial \Omega_\varepsilon}=0$) or Neumann ($\partial_n u_{|\partial\Omega_\varepsilon}=0$ where $\partial_n$ is the derivative normal to the boundary) conditions. We assume boundedness of $\Omega_\varepsilon$ and sufficient regularity of its boundary to ensure that the spectrum is discrete with finite multiplicities. Therefore, for each $\varepsilon>0$ we obtain a sequence
$$ \lambda_1(\varepsilon)\leq \lambda_2(\varepsilon)\leq \dots \to \infty$$
of eigenvalues, where each eigenvalue appears as often as the dimension of its eigenspace indicates. Corresponding eigenfunctions will generally be denoted $u_{k,\varepsilon}$.

\begin{main problem}
Describe the asymptotic behavior of the eigenvalues $\lambda_k(\varepsilon)$, of spectral invariants constructed from these, and of the eigenspaces as $\varepsilon\to 0$.
\end{main problem}

We now describe some other sorts of thin tubes, which will be referred to throughout the article.
\begin{itemize}
\item {\bf Manifolds and higher dimensions: } Instead of plane domains one may consider domains in $\R^n$ and more generally Riemannian manifolds, with the Laplace-Beltrami operator.\footnote{Instead of the Laplace-Beltrami operator, one may also consider other operators relating to the tube geometry, for example the Dirac operator.} Thus, for some $n\geq 2$ we start with an $(n-1)$ dimensional manifold $Y$ and two $n$-dimensional manifolds $X_l,X_r$. All manifolds are compact and Riemannian and may have (piecewise smooth) boundary.\footnote{for an exact definition of this, see for example \citet*{HasZel:QEBVE}; for our purposes one may assume that the data are such that the spaces defined below are smooth manifolds with boundary.} Also, we are given  isometries $j_l,\ j_r$ from $Y$ to subsets of the boundary of $X_l$, $X_r$. For $\varepsilon>0$ and a Riemannian manifold $(X,g)$ we denote by $\varepsilon X$ the Riemannian manifold $(X,\varepsilon^2 g)$, that is, all lengths are multiplied by $\varepsilon$. For an interval $I=[a,b]$ we then set
\begin{equation}
\label{eqn Omega3}
\Omega_\varepsilon =\left[\varepsilon X_l\, \cup\, (I\times \varepsilon Y)\, \cup\, \varepsilon X_r\right] / \{j_l,j_r\}
\end{equation}
({\em the thin tube manifold with cross section $Y$ and ends $X_l,X_r$})
where the union is disjoint and the quotient means identifying $a\times Y$ with  $j_l(Y)\subset\partial X_l$ using $j_l$ and $b\times Y$ with $j_r(Y)\subset\partial X_r$ using $j_r$.

Note that if $Y$ is an interval and $X_l,X_r\subset\R^2$ then this reduces to \eqref{eqn def Omega2} with $h_1$ and $h_2$ constant. It is straightforward to give a manifold setting where the cross-section has varying thickness.

Note also that this setting includes the $\varepsilon$-neighborhood of a straight-line segment in $\R^n$ (that is, points of distance at most $\varepsilon$ from the line), as well as the boundary of this neighborhood (sometimes referred to as sleeve). In the former case, $X_l,X_r$ are half balls and $Y$ is the $(n-1)$-dimensional unit ball, and in the latter case $X_l,X_r$ are half spheres and $Y$ is a sphere.

The Laplace(-Beltrami) operator on a space $X$, with boundary conditions given in the context, will be denoted by $\Delta_X$. We always take it positive, i.e.\ in $\R^n$ it is
$$ \Delta_{\R^n} = -\sum_{i=1}^n \frac{\partial^2}{\partial{x_i}^2}$$
\item {\bf Curved and twisted tubes: } Instead of modeling $\Omega_\varepsilon$ on a straight line, one may consider domains obtained as $\varepsilon$-neighborhood of an arbitrary simple curve. More generally, instead of a circular cross-section one may look at arbitrary cross-sections. The cross section may also rotate when progressing along the curve.
Such tubes may be considered with or without ends.

\item {\bf Tubes modeled on graphs: } An important generalization of tubes with ends is obtained by replacing the single line segment by a metric graph. Here, a graph is given by sets $V$ of vertices and $E$ of edges, where each edge connects two vertices (or one vertex with itself). Write $e\sim v$ if edge $e$ has vertex $v$ as an endpoint. A metric graph is a graph together with  a function $l:E\to(0,\infty)$. $l(e)$ is to be thought of as the length of $e$. Suppose we are given an $(n-1)$-dimensional manifold $Y_e$ for  each edge $e$ and an $n$-dimensional manifold $X_v$ for each vertex $v$, as well as isometries (gluing maps) of $Y_e$ with a subset of the boundary of $X_v$ for all $e\sim v$, without overlaps (i.e., the images of the gluing maps corresponding to a single vertex are disjoint). Then $\Omega_\varepsilon$ is defined by gluing cylinders of length $l(e)$ and cross section $\varepsilon Y_e$  to
$\varepsilon X_v$ for all pairs $e\sim v$.

A special case is the $\varepsilon$-neighborhood of a graph embedded in $\R^n$ with straight edges, and the boundary of this neighborhood. Of course, one may also consider curved and twisted edges.

The additional combinatorial structure of the graph suggests a refinement of the Main Problem in this context:

\begin{Qgraph problem}
How is the limiting behavior of the spectral data on $\Omega_\varepsilon$ related to spectral data of some operator on the metric graph itself?
\end{Qgraph problem}
 This question naturally leads to quantum graphs, to be introduced below.
\smallskip

\item
{\bf Noncompact tubes: }So far, all tubes were compact. Instead, one may consider tubes modeled on the real line or half line, or more generally on graphs having some edges of infinite length (which then have only one endpoint). Such spaces arise in their own right and also as certain rescaled limit objects of families of compact tubes, see equation \eqref{eqn def Xinfty}. Apart from this, however, we focus mostly on compact tubes in this article.
\nocite{Pav:SGMOE}
\end{itemize}

The goal of this article is to explain some common ideas and methods used in different areas where thin tubes occur. We also explain a very limited selection of  mostly  recent results, mainly from mathematical physics and spectral geometry.

Thin tubes as defined above are a special class of singular degenerations. There are many other kinds of singular degenerations. For example, one may keep the vertex neighborhoods fixed while only shrinking the tubes between them (this is sometimes called a handlebody), or one may remove a small ball inside some fixed domain. There is a large body of literature on such problems. We just mention the book by \citet*{MazNazPla:ATERSGG}. See also the discussion in \ref{subsec pseudo}.

In addition to the references given below, we would like to mention the collection by \cite{BerCarFul:QGTA} which has many articles on related questions.

The author acknowledges the hospitality and support of the Isaac Newton Institute in Cambridge during the program 'Analysis on Graphs and its Applications' in 2007, where part of this work originated.

\subsection{Rescaling to long cylinders}
For a proper understanding of many tube problems it is important to consider tubes rescaled to fixed width. Denote
\begin{equation}
\label{eqn def rescaling}
N=\varepsilon^{-1},\quad
X^N = \varepsilon^{-1}\Omega_\varepsilon.
\end{equation}
That is, all lengths in $\Omega_\varepsilon$ are stretched by the factor $\varepsilon^{-1}$. This gives a tube of fixed width and with fixed ends (independent of $\varepsilon$) and length on the order of $N=\varepsilon^{-1}$ (assuming $\Omega_\varepsilon$ is given as in \eqref{eqn def Omega2} or \eqref{eqn Omega3}, and similarly for tubes modeled on a graph).
The eigenvalues and eigenfunctions of the Laplacian on $X^N$ are\footnote{\label{fn rescale}We will denote rescaled quantities and variables with a bar, for example $\lambdabar$, $\ubar$ for eigenvalue and eigenfunction, and $\xbar,\ybar$ for coordinates on $X^N$. }
\begin{equation}
\label{eqn eigenvalue rescaling}
\lambdabar_k(N) = \varepsilon^2 \lambda_k(\varepsilon),\quad \ubar_{k,N}(\xbar,\ybar) = u_{k,\varepsilon}(\varepsilon\xbar,\varepsilon\ybar).
\end{equation}
Note that the limit $\varepsilon\to 0$ corresponds to $N\to\infty$, and this yields a limit object $X^\infty$ different from $\Omega_0$. While $\Omega_0$ is just the one-dimensional skeleton of the tube, $X^\infty$ has the same dimension as the tube and also knows about the shape of the cross sections and of the ends.
In this sense, $X^\infty$ retains more information about the tubes than $\Omega_0$. On the other hand, in the case of graphs, $\Omega_0$ knows about the graph (connectivity) structure while $X^\infty$ does not.
The precise definition of $X^\infty$ will be given in \eqref{eqn def Xinfty}.

In fact, in part of the literature it is common (and more natural) to state problems and results in the 'long cylinder' setting, especially in the global analysis literature.

In this article I will state most results in terms of thin tubes, so that they can be compared more easily.

\section{Origins} \label{sec origins}
\subsection{Mathematical Physics}
There are a number of physical systems which are quasi one-dimensional in the sense that their dimensions in one direction are much larger than their dimensions in transversal directions, or more generally that several such one-dimensional subsystems are coupled in a graph-like structure. The quantum mechanics of such systems is then most easily modeled by the underlying graph (i.e.\ neglecting all transversal directions) with a Schrödinger operator on the edges coupled by transition (boundary) conditions at the vertices (a so-called quantum graph). A better, still simplifying model, is to consider thin tubes around the graph edges and a Schrödinger operator, for example the Laplacian on such tubes. Examples of such systems are: carbohydrate molecules, where the graph structure is given by the nuclei and their bonds -- which are assumed fixed --, and the goal is to analyze the behavior of the electrons; this in fact seems to be the earliest occurrence of such a model, see \citet*{RueSch:FENMCSI}. Another, more recent kind of such quasi one-dimensional systems that need to be analyzed quantum-mechanically are so-called nanotubes or quantum wires, or highly integrated circuits. See the nice survey by \citet*{Kuc:GMWTS} for more background and many references, and the classical article by \citet*{ExnSeb:ESMCOT} which sparked a lot of activity in this field.

Confining a quantum-mechanical particle to a small neighborhood of a graph-like structure may be achieved also by other means than by a thin tube: For example, one may consider the Schrödinger operator on $\R^n$ with a potential $a^2 V$ added, where $a$ is a large parameter and $V$ is non-negative, and vanishes precisely on the graph. This leads to closely related questions and results. See for example \cite{SmiSol:QGLNPW}, where a transversal harmonic oscillator is added to the quantum graph operator, and \cite{FroHer:RHCCQM} (in the case of submanifolds instead of graphs).

\subsection{Global analysis}
Global analysis studies how analytic objects associated with certain partial differential operators on a manifold (e.g.\ spectrum, dimension of kernel) are related with global, e.g.\ topological, properties of the manifold. This includes for instance Hodge theory (where the cohomology of a compact manifold is calculated in terms of the kernel of the Laplace-Beltrami operator on differential forms) and index theory (where the difference of the  dimensions of the solution spaces of a system of PDEs on a compact manifold and of its adjoint is computed by topological means). A typical procedure is this: one starts with a smooth compact manifold $X$, then chooses a Riemannian metric $g$ on $X$, this defines the Laplace-Beltrami operator $\Delta_X$ (or, in the presence of additional structure, Dirac operators), then one forms certain expressions from the spectrum of $\Delta_X$, and finally these are shown to be independent of the choice of metric, and hence  (dif\-fe\-ren\-tial-)\-to\-po\-lo\-gi\-cal invariants. This happens in Hodge theory, and another instance of this is the analytic torsion, which is a combination of determinants of the Laplacians (introduced below in \ref{subsubsec global}) on $k$-forms, for $k=0,\dots,\dim X$. The determinant of the Laplacian itself is  only a {\em spectral invariant}, i.e.\ it depends on the eigenvalues (and hence on the metric).
There are also other spectral invariants (for example, the $\eta$-invariant for Dirac-operators) which appear in global analysis contexts.

It is then a basic problem to understand how these topological invariants change under modifications of the manifold: For example, if one cuts the manifold into two pieces by removing a hypersurface $Y$ then one would like to know how the invariant of the whole manifold relates to the invariants of the pieces (one seeks a {\em gluing formula}).

The same question can be asked for spectral invariants. This is interesting in its own right, but since the topological invariants are defined in terms of these,
an answer to the spectral question will also answer the topological question.

Surprisingly, a useful technique for obtaining such gluing formulas is to consider a family of Riemannian manifolds $X^N$ obtained by gluing a cylinder\footnote{to obtain a smooth metric, one has to assume that the metric on $X$ is of product type, i.e.\ $dx^2+g_Y$, in a neighborhood of $Y$} $[-N,N]\times Y$ into $X$ in place of $Y$, and studying the asymptotic behavior of the spectral invariants of $X^N$
as $N\to\infty$, in terms of those of $Y$ and of the ends $X_l$, $X_r$. An example of this will be given in \ref{subsubsec global}. As manifold, $X^N$ is the same as (more precisely, diffeomorphic to)  $X$. It is only the metric which is different. Since the topological invariants are independent of the metric, hence of $N$, this asymptotic behavior determines the behavior also for $N=0$, thus providing the desired gluing formula.

\subsection{Spectral geometry}
Spectral geometry studies how the eigenvalues and eigenfunctions of the Laplacian on a space depend on the geometry of the space, and which properties they have independent of the geometry. In this context, thin tubes sometimes serve as a special class of spaces for which questions that are very hard in general may be answered by asymptotics methods. We give two examples of this:

\medskip

\noindent{\bf The nodal line conjecture}
\\
We consider the Dirichlet Laplacian on a bounded plane domain $\Omega$. It is well-known that the first (smallest) eigenvalue is simple, that the corresponding eigenfunction does not change sign in $\Omega$ and that consequently any eigenfunction $u$ for the second smallest eigenvalue must have a non-empty {\em nodal line} $\calN:=u^{-1}(0)$. Also Courant's nodal domain theorem implies that $\calN$ divides $\Omega$ into two connected components, and this implies that it is a smooth curve.
\citet*{Pay:TCFMEP} conjectured that $\calN$ must touch the boundary of $\Omega$, that is, it cannot be a closed curve in the interior of $\Omega$.\footnote{The state of the nodal line conjecture is as follows: \citet*{Mel:NLSELR} proved it for convex plane domains. \citet*{HofHofNad:NLSELRCC} showed that it cannot hold if one allows the domain to be multiply connected. \citet*{Jer:FNSCD} proved it for convex domains in $\R^n$ with diameter one and small inradius. An interesting higher dimensional result was also obtained by \citet*{FreKre:LNSTCT}. An open problem is whether the nodal line conjecture is true for plane simply connected domains.}

\citet*{Jer:FNLCPD} approached the problem in the case of convex domains $\Omega$ as follows: By scale invariance, one may assume that the diameter of $\Omega$ (the largest distance between two points in $\Omega$) is one. Let $\varepsilon$ be the inradius of $\Omega$, that is the maximum radius of a ball contained in $\Omega$. It is not  hard to show that $\Omega$ is of the form \eqref{eqn def Omega1}, where $h_1,h_2$ are uniformly bounded with $\max (h_2-h_1)$ of order one.

Therefore, if the inradius is small then this is a case of thin tubes, and therefore some of the techniques described below are applicable in principle. However, a major complication arises from the fact that the first and second derivatives of $h_1,h_2$ cannot be controlled beyond what follows from convexity. In spite of this, Jerison manages to prove the nodal line conjecture for convex domains with sufficiently small inradius.
See \ref{subsubsec variable low reg} for more on this.

\medskip

\noindent{\bf Multiplicity of the second eigenvalue}
\\
A basic question in spectral geometry is: Given  a compact manifold $M$ without boundary and an increasing sequence of non-negative real numbers, is there a Riemannian metric on $M$ such that the given sequence is precisely the eigenvalue sequence of the Laplace-Beltrami operator? This is probably intractable, so one may simplify the question by just prescribing a finite number, say the first $N$ eigenvalues (always counted with multiplicities).

Assuming $M$ is connected and without boundary, a first condition on the sequence is that the first eigenvalue be zero and have multiplicity one.
In two dimensions there is another well known obstruction: The multiplicity of the higher eigenvalues is bounded in terms of the Euler characteristic of $M$, see \citet*{Che:ENS}.

In dimensions $n\geq 3$, \citet*{Col:MPVPNNL} showed that there is no such obstruction for the multiplicities of the second eigenvalue\footnote{In fact, \citet*{Col:CLDPFSD} shows the stronger result that any sequence of $N$ numbers $0=\lambda_1<\lambda_2\leq\dots\leq\lambda_N$ can occur.}. The idea of the proof is as follows: Given $N\in\N$, embed a complete graph with $N$ vertices in $M$ (here one needs that $\dim M\geq 3$, if $N$ is large). On a neighborhood of this graph in $M$, define the metric to be of the thin tube type $\Omega_\varepsilon$ (with all edges of length one and disks as cross sections), and extend it outside of this neighborhood so that it is small there in a suitable sense. Then the eigenvalues on $M$ will be close to those on $\Omega_\varepsilon$, with Neumann boundary conditions at $\partial\Omega_\varepsilon$, and  these will be close to those of a suitable operator on the graph itself if $\varepsilon$ is small (in the latter step the methods discussed below are used; see below for the definition of the operator on the graph). A simple calculation shows that the second eigenvalue of the operator on the graph has multiplicity $N-1$. Hence the second to $N$th eigenvalue of $M$ must lie very close together. A topological perturbation argument then shows that one may change the metric on $M$ to obtain one where the multiplicity of the second eigenvalue equals $N-1$.

\section{Results}
We will explain a number of results and try to motivate why they should be true. The selection is restricted by space as well as the personal taste and the limitations of the author. I apologize to anyone (that is, almost everyone) whose result is not included.

We will focus on results directly related to spectral quantities. In some cases, there are actually stronger statements in the papers concerning the resolvent, which imply the spectral results. See the explanation in \ref{subsec resI}. We emphasize results where Dirichlet (or other non-Neumann) boundary conditions are imposed, since this case has seen the greatest progress recently. Ideas of proofs of some of the results are given in Section \ref{sec methods}.

After introducing the basic principle that guides what one should expect for thin tube problems we consider results where the limit object is an interval or a curve and then the case of graphs. Clearly intervals are a special case of graphs. However, there are many results that are only natural for intervals, or have so far only been formulated for these, or are better (for example, sharper) than their graph analogues.

\subsection{The basic principle} \label{subsec basic principle}
The basic principle in thin tube problems is that for the lower part of the spectrum the first transversal mode dominates the behavior of eigenvalues and eigenfunctions.
To see this, we first look at the simple example of tubes modeled on an interval $I$ and with constant cross section $Y$.
Let
\begin{equation}
\label{eqn def nu}
\nu_0\leq \nu_1\leq \dots \text{ the eigenvalues of }\Delta_Y,\text{ with eigenfunctions }\varphi_0,\varphi_1,\dots,
\end{equation}
with the boundary conditions given at $\partial Y$.
Then  the eigenfunctions and eigenvalues on  $\Omega_\varepsilon=I\times \varepsilon Y$ are\footnote{A remark on notation: If $Y$ is a Riemannian manifold, instead of a domain in $\R^{n-1}$, then $y/\varepsilon$ does not really make sense. We use it anyway since it is suggestive. The precise meaning is: $y$ is a point of $Y$ considered with the rescaled metric $\varepsilon^2 g_Y$, and $y/\varepsilon$ is the same point in $Y$, where we now use the metric $g_Y$.}
\begin{equation}
\label{eqn product case}
u_\varepsilon(x,y) = \sin (a_l x)\, \varphi_m(\frac y\varepsilon),\quad
\lambda(\varepsilon) = \varepsilon^{-2}\nu_m + a_l^2,\ a_l=\frac{\pi l}{|I|},
\end{equation}
for $l\in\N$, $m\in\N_0$ (for Dirichlet conditions at the ends of $I$; replacing $\sin$ by $\cos$ and taking $l\in\N_0$ gives Neumann conditions). So for any $k\in\N$, the $k$th eigenvalue is $\varepsilon^{-2}\nu_0 + a_k^2$ as soon as $\varepsilon$ is sufficiently small.
This leads us to expect for general tubes that
\begin{itemize}
\item
for the (rough) behavior of $\lambda_k(\varepsilon)$ for {\em fixed} $k$ and $\varepsilon\to 0$ only the first transversal mode ($m=0$) should be relevant, and
\item
higher transversal modes ($m>0$) will be relevant in the behavior of $\lambda_k(\varepsilon)$ if, simultaneously with $\varepsilon\to0$, one lets $k\to\infty$ on the order of $\varepsilon^{-1}$.
\end{itemize}
In the case of a plane domain with variable cross section, \eqref{eqn def Omega1}, the first eigenvalue of the transversal operator at $x$, $-\partial^2/\partial \ybar^2$ on the interval $[h_1(x),h_2(x)]$ with Dirichlet conditions, is $\pi^2/h(x)^2$ with $h(x)=h_2(x)-h_1(x)$, and correspondingly one expects the behavior of low eigenvalues to be dominated by those of the one-dimensional Schrödinger operator
\begin{equation}
\label{eqn schroedinger}
P_0 = -\frac{d^2}{d x^2} + \varepsilon^{-2} \frac{\pi^2}{h(x)^2},\quad x\in I,
\end{equation}
(with boundary conditions at $\partial I$ corresponding to those of the original problem) and the eigenfunctions to be approximately
\begin{equation}
\label{eqn eigenfcn approx}
u_\varepsilon(x,y)\approx u_0(x) \sin \frac{\pi}{\varepsilon h(x)} (y-h_1(x)),
\end{equation}
for eigenfunctions $u_0$ of $P_0$.

There are many results that make these expectations precise in one way or another. They differ in the generality of the setup (constant or variable width; tubes with or without ends; interval, curve or graph as limit object; regularity assumptions etc.) as well as in the precision that is aimed at. The desired precision is motivated by the origin of the problem: In some problems interesting behavior occurs already in the leading order asymptotics, while for some one has to go to second or third order (in which case it may be natural to consider asymptotics to arbitrary order). Having higher order asymptotics may also be useful since in some cases it simplifies the proof of the validity of the lower orders. In reverse, this means that in irregular problems (where there can be no higher order asymptotics) it is technically harder to prove validity of even leading order asymptotics. An example of this is \citep*{GriJer:AFNLCD,GriJer:SFECPD}, see Section \ref{subsubsec variable low reg}.

\subsection{Intervals and curves as limits: Tubes without ends}
As will be seen when considering graphs (see \ref{subsubsec graphs Neumann}) the case of Neumann conditions is rather uninteresting if $\Omega_0$ is a single curve. Therefore, it is not surprising that all the results mentioned here assume Dirichlet boundary conditions. There is one exception concerning nodal lines in \ref{subsubsec variable low reg}.
\subsubsection{Domains of variable thickness, I} \label{subsubsec varthick I}
 \citet*{FriSol:SDLNS,FriSol:SDLNSII} consider a domain \eqref{eqn def Omega1} with $h_1=0$, and they assume that $h=h_2$ has a unique maximum on $I$, say at the interior point $x=0$, and that near $x=0$ it has $p$th order behavior (normalized)
\begin{equation}
\label{eqn pth order}
h(x) = \begin{cases}
1-c_+  x^p + O(x^{p+1}),& \quad x>0,\\
1-c_-  |x|^p + O(|x|^{p+1}),& \quad x<0,
\end{cases}
\end{equation}
while it is $C^1$ on $I\setminus\{0\}$. Here $c_\pm >0$ and $p\geq 1$. For example, for $p=1$ a triangle is allowed. One considers Dirichlet boundary conditions at the upper and lower boundary.

What should we expect in this setting? According to the basic principle of tube problems one expects the behavior of eigenvalues to be dominated by the Schrödinger operator \eqref{eqn schroedinger}. This is, in terms of its dependence on the small parameter $\varepsilon$, a semi-classical operator (replace $\varepsilon$ by $\hbar$), and for these it is standard that eigenfunctions concentrate strongly near the minimum of the potential $V(x)=\frac{\pi^2}{h^2}$, i.e. near $x=0$. So the values of $V$ away from $x=0$ should have 'small' influence on the behavior of eigenvalues and eigenfunctions. Therefore, we replace $V(x)=\pi^2\left(1 + 2c_\pm |x|^p + O(|x|^{p+1}\right)$ (from \eqref{eqn pth order}) on $I$ by $\pi^2 + 2\pi^2 c_\pm |x|^p$ on $\R$. The first summand gives a contribution of $\varepsilon^{-2}\pi^2$ to the eigenvalue, and the second a contribution equal to an eigenvalue of the Schrödinger operator $-\frac{d^2}{dx^2} + \varepsilon^{-2}C_\pm |x|^p$, where $C_\pm=2\pi^2 c_\pm$. Rescaling $x=\varepsilon^{\alpha}\xi$, $\alpha=\frac2{p+2}$, shows that these eigenvalues equal $\varepsilon^{-2\alpha}\mu_k$ where
\begin{equation}
\label{eqn def muk}
\mu_k\text{ are the eigenvalues of } -\frac{d^2}{d\xi^2}+C_\pm |\xi|^p\text{ on }\R.
\end{equation}
Summarizing, one should expect
\begin{equation}
\label{eqn result FriSol}
\lambda_k(\varepsilon) = \varepsilon^{-2}\pi^2 + \varepsilon^{-2\alpha} \mu_k + o(\varepsilon^{-2\alpha}),\quad \alpha=\frac2{p+2}.
\end{equation}
This is precisely the result of \citet*{FriSol:SDLNS} (in the case of Dirichlet conditions at the left and right ends) and \citet*{FriSol:SDLNSII} (in the case of Neumann conditions at the left and right ends, and in a similar form with the interval $I$ replaced by $\R$). Their main tool is a resolvent analysis, see \ref{subsec resI}. They also prove an approximation \eqref{eqn eigenfcn approx} of the eigenfunctions (in $L^2$-norm), where $u_0$ can be replaced by the appropriately rescaled eigenfunctions of the operator in \eqref{eqn def muk}.

\subsubsection{Domains of variable thickness, II: Low regularity}\label{subsubsec variable low reg}
\citet*{Jer:DFNLCD} and \citet*{GriJer:AFNLCD,GriJer:SFECPD} consider domains \eqref{eqn def Omega1} but only assume
\begin{itemize}
\item convexity, i.e.\ $h_1$ is convex and $h_2$ is concave,
\item uniform bounds on $h_1,h_2$; we normalize\footnote{it is not hard to show that the inradius is $\varepsilon+o(\varepsilon)$ then} these so
$$ \max_{x\in I} h(x) = 1,\quad h=h_2-h_1.$$
\end{itemize}
There are no bounds on the derivatives of $h$, it is not even assumed to be differentiable.\footnote{This complicates the analysis considerably since when applying $\Delta$ to \eqref{eqn eigenfcn approx} then one obtains, in addition to the terms yielding $P_0 u_0$, terms involving first and second derivatives of $h$ and $h_1$, and the point is to control these.} Put differently, one does not fix $h_1,h_2$ and then lets $\varepsilon\to 0$, but rather one considers, for any small $\varepsilon$, all possible $h_1,h_2$. Note that in this class the distinction 'with end' - 'without ends' is immaterial since an end (in the convex case, say) corresponds to $h_{1/2}$ having derivative on the order of $\varepsilon^{-1}$ in an $\varepsilon$-neighborhood of $\partial I$ (plus a trivial $1+C\varepsilon$ rescaling of length).
One assumes Dirichlet boundary conditions.
\medskip

{\em Eigenvalues: } Can one expect a two-term asymptotics for the eigenvalues? Since in this setting one finds domains of the form \eqref{eqn pth order} with arbitrary $p$, one may expect from \eqref{eqn result FriSol} that the answer is no.

Surprisingly, one may derive from $h_1,h_2$ a geometric quantity which almost yields a second term for $\lambda_{k}(\varepsilon)$:\footnote{In the cited papers only $k=1,2$ are considered, but at least for the eigenvalues the results extend to arbitrary $k$ easily.}
\citet*{Jer:DFNLCD} proves that there are absolute constants $c,C>0$ such that
\begin{equation}
\label{eqn Jer result1}
\varepsilon^{-2}\pi^2 + \frac c{\ell^2} \leq \lambda_k(\varepsilon) \leq
\varepsilon^{-2}\pi^2 + \frac C{\ell^2}
\end{equation}
where $\ell$ is the biggest number such that
\begin{equation}
\label{eqn def l}
h \geq 1 - \frac{\varepsilon^2}{\ell^2}\text{ on an interval }I_\ell\subset I\text{ of length }\ell.
\end{equation}
For example, if $I=[0,1]$ then $\ell=1,\ I'=I$ for a rectangle ($h_1,\ h_2$ constant) while $\ell=\varepsilon^{2/3},\ I=[0,\ell]$ for a triangle ($h_1\equiv 0$, $h_2(x)=1-x$), and it can be shown that $\ell$ always lies between these values. Compare also with \eqref{eqn result FriSol} where $\ell\approx\varepsilon^{\alpha}$. $\ell$ should be regarded as a characteristic length scale for $\Omega_\varepsilon$.
\medskip

{\em Eigenfunctions:}
To which precision does the approximation \eqref{eqn eigenfcn approx} hold? More generally, to what precision can qualitative properties of $u_\varepsilon$ be deduced from a knowledge of $h$ alone?\footnote{Note that different domains can have the same height function $h$. This question may be understood in two ways: First, one uses only 'elementary' (geometric) information about $h$. Second, one allows to use transcendental information, for example the eigenfunctions $u_0$ of the Schrödinger operator \eqref{eqn schroedinger}. We consider both questions.} Normalize $u_\varepsilon$ by assuming $\max |u_\varepsilon|=1$, \citet*{Jer:DFNLCD} shows that $u_\varepsilon$ 'lives' near the interval $I_\ell$ in \eqref{eqn def l} in the sense that it decays like $e^{-d(x)/\ell}$ in terms of the distance $d(x)$ from $I_\ell$. \citet*{GriJer:SFECPD} study the first eigenfunction ($k=1$) in detail. The approximation \eqref{eqn eigenfcn approx} is shown to hold on and near $I_\ell$ with an error $O(\varepsilon/\ell)$, and, if $\locmax f$ denotes the point where a function $f$ assumes its maximum (this is unique for the first eigenfunctions $u_1,\varphi_1$ in the case of convex domains) then
\begin{equation}
\label{eqn locmax}
\locmax u_{1,\varepsilon} = (\locmax \varphi_1,0) + O(\varepsilon).
\end{equation}
\citet*{GriJer:AFNLCD} study the second eigenfunction, in particular the nodal line, in detail. First, the approximation \eqref{eqn eigenfcn approx} is quantified to hold with an error $O(\varepsilon^2/\ell)$ near the nodal line (eqn.\ (15) in loc.cit.), and this is used to show that the nodal line is contained in a vertical strip of width $O(\varepsilon^2)$. All these estimates are sharp in order of magnitude. See Subsection \ref{subsubsec sharp loc max} for more on this.
\citet*{Jer:LFNLNP} studies the first nodal line in the Neumann problem from a similar perspective.

\subsubsection{Curved and twisted tubes}
Given a smooth embedded curve $\gamma:I\to\R^n$ parametrized by arc length, a bounded open set $Y\subset\R^{n-1}$ and a family $(\Phi_x)_{x\in I}$ of linear isometries $\Phi_x:\R^{n-1} \to N_{\gamma(x)}\gamma$ (the normal space to $\gamma$ in the point $\gamma(x)$) depending smoothly on $x$ we call the set
$$ \Omega_\varepsilon = \{\gamma(x) + \varepsilon \Phi_x(y):\, x\in I,\ y\in Y\}$$
a {\em thin curved tube with cross section $Y$}. It is embedded for small $\varepsilon$. We call it {\em non-twisted} if $\Phi$ can be chosen such that for each $y$ the curve $x\mapsto \Phi_x(y)$ cuts each normal slice $\gamma(x)+\Ran\Phi_x$ orthogonally. It is easy to see that in this case the (Euclidean) Riemannian metric on $\Omega_\varepsilon$, expressed in the coordinates $(x,y)\in I\times Y$, is diagonal of the form\footnote{Write $x_0=x$, $(x_1,\dots,x_{n-1})=y$. Then $g_{ij}=\partial_{x_i}\Psi\cdot\partial_{x_j}\Psi$ where $\Psi(x,y)=\gamma(x)+\varepsilon\Phi_x(y)$.  Since $\Phi_x$ is an orthogonal linear map for each $x$, we have $g_{ii}=\varepsilon^2$ for $i\geq 1$ and $g_{ij}=0$ for $i,j\geq 1$, $i\neq j$. By construction $\gamma'(x)\perp \partial_{y_j}\Phi_x(y)$  and by the non-twisting assumption $\partial_x\Phi_x(y)\perp \partial_{y_j}\Phi_x(y)$ for all $x,y$. Therefore $g_{0j} = (\gamma' + \varepsilon \partial_x\Phi)\cdot \partial_{y_j}\Phi=0$. Finally, by non-twisting, $\gamma'$ and $\partial_x\Phi$ are parallel, so  $g_{00}= \|\gamma'+\varepsilon \partial_x\Phi\|^2 = (1+\varepsilon\gamma'\cdot\partial_x\Phi)^2$, and differentiating $\gamma'(x)\cdot\Phi_x(y)\equiv 0$ in $x$ yields $\gamma'\partial_x\Phi + \gamma''\Phi= 0$ and hence the given expression for $g_{00}=a^2$.}
\begin{equation}
\label{eqn tube metric}
g= \diag(a^2,\varepsilon^2,\dots,\varepsilon^2),\quad a(x,y) = 1-\varepsilon \gamma''(x)\cdot\Phi_x(y).
\end{equation}
For example, any curved tube with $Y$ a disk centered at zero is non-twisted; also, given 'generic' $\gamma$ and $Y$, there is a non-twisted curved tube, unique if fixed at one end. See \citet*{FreKre:LNSTCT} and the references given there.

What should one expect for the eigenvalues of a curved, non-twisted tube with constant cross section, with Dirichlet boundary conditions?
The Laplacian on $\Omega_\varepsilon$ expressed in the coordinates $(x,y)$ is the variable coefficient operator\footnote{Use the general formula $\Delta_g = -\frac1{\sqrt{\det g}}\partial_{x_i} g^{ij}\sqrt{\det g} \partial_{x_j}$, where $(g^{ij})$ is the inverse matrix to $g$, and the volume element $\sqrt{\det g}dx_0\cdots dx_{n-1}$.}
$$\Delta_g = -a^{-1} \partial_x a^{-1}\partial_x - \varepsilon^{-2} a^{-1}\sum_{j=1}^{n-1} \partial_{y_j}a\partial_{y_j}$$
which is symmetric on its Dirichlet domain in $L^2(I\times Y,a\,dxdy)$. Then $P=a^{1/2}\Delta_g a^{-1/2}$ is a unitarily equivalent operator in $L^2(I\times Y,dxdy)$. Symmetry implies that it has no first order terms, and indeed from \eqref{eqn tube metric} one sees that
\begin{equation}
\label{eqn curved tube operator}
P=-a^{-2}\partial_x^2 + \varepsilon^{-2}\Delta_Y + V_\varepsilon
\end{equation}
where $\Delta_Y$ is the standard (Euclidean) Laplacian on $Y$ and $V_\varepsilon(x,y)=-\frac14 \kappa(x)^2 + O(\varepsilon)$
is smooth in $\varepsilon,x,y$. Here $\kappa(x)=|\gamma''(x)|$ is the curvature. Since $a=1+O(\varepsilon)$ this is a $O(\varepsilon)$ perturbation of $(-\frac{d^2}{dx^2} - \frac14 \kappa^2) + \varepsilon^{-2}\Delta_Y$, which is a direct sum operator on $I\times Y$. Therefore, one expects the eigenvalues of $P$, and hence of $\Delta$ on $\Omega_\varepsilon$, to be
\begin{equation}
\label{eqn ev on curved tube}
\lambda_k(\varepsilon) = \varepsilon^{-2}\nu_0 + \mu_k + O(\varepsilon),\quad \mu_k=\text{ eigenvalues of }-\frac{d^2}{dx^2}-\frac {\kappa(x)^2}4.
\end{equation}
This is precisely the result of \citet*{FreKre:LNSTCT} (see also the nice survey by \cite{Kre:TVBQW}). They also consider slightly more general operators of the form \eqref{eqn curved tube operator}. In addition they show (assuming $\partial Y$ is smooth) that the eigenfunctions are $O(\varepsilon)$-close to the product type eigenfunctions $\psi_k(x) \varphi_0(y/\varepsilon)$ (with $\psi_k$ the eigenfunction corresponding to $\mu_k$), uniformly (!) with error linearly decaying at the boundary, as well as certain uniform derivative estimates. They use this to show that the nodal set (zero set) of the $k$-th eigenfunction on $\Omega_\varepsilon$ must lie within $O(\varepsilon)$ of the nodal set of $\psi_k$ and that it must touch the boundary for small $\varepsilon$.

Curved and twisted tubes were studied in many papers, for example by \citet*{ExnSeb:BSCQW}, \citet*{DucExn:CIBSQWTTD}, \citet*{EkhKovKre:HITW}, \citet*{BouMasTra:CTEODW}, mostly as waveguides.  See also \citet*{GeeKel:ESCWST} for a (non-rigorous) WKB-analysis of curved waveguides and tubes with variable cross section.

\subsection{Intervals and curves as limits: Tubes with ends}
What difference does an end make? Let us see what to expect, by adding an $(\alpha\varepsilon)\times\varepsilon$ rectangle as an end to a $1\times\varepsilon$ rectangle, for a fixed $\alpha>0$. This yields a $(1+\alpha\varepsilon)\times\varepsilon$ rectangle, so the eigenvalues change from $\varepsilon^{-2}\pi^2 m^2 + \pi^2 l^2$ to
$\varepsilon^{-2}\pi^2 m^2 + \pi^2 l^2(1+\alpha\varepsilon)^{-2} = \varepsilon^{-2}\pi^2 m^2 + \pi^2 l^2 - \varepsilon\,2\alpha\pi^2 l^2+O(\varepsilon^2)$, that is, the end makes a contribution of order $\varepsilon^1$, which is rather small and hence negligible for many problems.

However,  ends do make an important difference in a number of situations:
\begin{itemize}
\item
In some spectral geometry problems it is precisely the $\varepsilon^1$ term of the eigenvalue asymptotics which matters, so the ends matter. See \ref{subsubsec sharp loc max}.
\item
In the global analysis context the influence of the ends is actually the main focus (see the explanations in Section \ref{sec origins}). It is also stronger since the spectral invariants involve all eigenvalues simultaneously.
\item
As we will see in \ref{subsubsec graphs general}, for graphs the precise form of the ends may influence even the $\varepsilon^{-2}$ and $\varepsilon^0$ term in the eigenvalue asymptotics.
\end{itemize}

\subsubsection{Rectangles with an end}\label{subsubsec sharp loc max}
\begin{figure}
\includegraphics[scale=0.8]{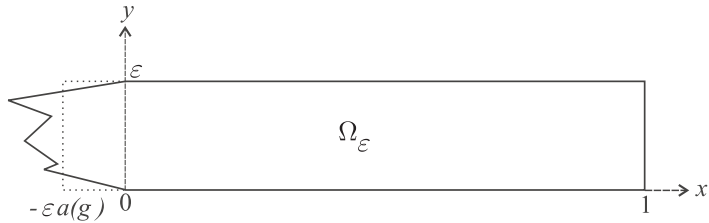}
\caption{A rectangle with an end}
\label{figure rectangle with end}
\end{figure}

\citet*{GriJer:AEPD} consider the Dirichlet problem on rectangular domains with one end added: Fix a Lipschitz function $g:[0,1]\to [0,\infty)$ and set
\begin{equation}
\label{eqn def phi GriJer}
\Omega_\varepsilon = \{(x,y)\in\R^2;\,y=\varepsilon\ybar,\ 0<\ybar<1,\ -\varepsilon g(\ybar)<x<1\}.
\end{equation}
See Figure \ref{figure rectangle with end}.
For the precise description of the eigenvalues and eigenfunctions, one needs to consider the rescaling \eqref{eqn def rescaling} and its limit as $N\to\infty$, defined as
\begin{equation}
\label{eqn rescaled limit GriJer}
X^\infty = \{(\xbar,\ybar)\in\R^2;\, 0<\ybar<1,\ -g(\ybar)<\xbar\}.
\end{equation}
This is a half infinite strip with an 'end' added at the left side.
The Laplacian on $X^\infty$ has continuous spectrum $[\pi^2,\infty)$. We need some information\footnote{It is also important that there are no $L^2$ eigenvalues $\leq\pi^2$ (see the explanations in \ref{subsec graphs}). This is true (and easy to see, using the variational principle) in this setting, because the end is contained in the doubly infinite strip $\R\times (0,1)$.} about the spectral value $\pi^2$:
\begin{quote}
There is a unique function $U$ on $\Xinfty$ satisfying $\Delta U=\pi^2 U$ in $\Xinfty$, $U=0$ at $\partial\Xinfty$ and $U(\xbar,\ybar)= \xbar\sin (\pi \ybar) + O(1)$.
\\
Furthermore, there is a unique number $a=a(g)\geq 0$ such that
\begin{equation}
\label{eqn U asymptotics}
U(\xbar,\ybar) = (\xbar+a(g))\sin(\pi \ybar) + O(e^{-\xbar}),\quad \xbar\to \infty.
\end{equation}
\end{quote}
$a(g)$ is related to the scattering phase, see footnote \ref{fn scatt phase}.
\eqref{eqn U asymptotics} may be interpreted as saying that, spectrally, $\Xinfty$ is similar to the half infinite strip $(-a(g),\infty)\times (0,1)$ with Dirichlet boundary conditions. Therefore, one may expect that the spectral data for $X^N$ should be similar to that for the rectangle $(-a(g),N)\times (0,1)$. This is in fact true to rather high precision. \citet*{GriJer:AEPD} prove (scaling back to $\Omega_\varepsilon$)
\begin{align}
\label{eqn GriJer lambdak}
\lambda_k(\varepsilon) &= \varepsilon^{-2}\pi^2 + \frac{k^2\pi^2}{(1+\varepsilon a)^2} + O(\varepsilon^3)\\
\label{eqn GriJer uk}
u_{k,\varepsilon}(x,y) & = \sin (k\pi\frac{x+\varepsilon a}{1+\varepsilon a})\,\sin (\varepsilon^{-1}\pi y) + O(\varepsilon^3) \quad\text{ for }x>3\varepsilon\log\varepsilon^{-1}.
\end{align}
for the suitably normalized $k$th eigenfunction, where the estimates are uniform (for fixed $k$) with all $\partial_{\xbar},\partial_{\ybar}$ derivatives.

Here is a question where a precision of estimates that 'detects' $a$ is needed. Consider the problem whether the estimate \eqref{eqn locmax} is sharp in order of magnitude, i.e.\ whether $O(\varepsilon)$ may be replaced by $O(\varepsilon^r)$ for some $r>1$. To prove a negative answer requires constructing an example which is far enough from product type but still can be analyzed sufficiently explicitly. This can be answered using \eqref{eqn GriJer uk}, and it is the quantity $a(g)$ which is essential here: \eqref{eqn GriJer uk} implies that $\locmax u_1= \frac12 - \varepsilon\frac {a(g)}2 + O(\varepsilon^2)$. So if one finds two functions $g_1,g_2$ inducing the same height function $h$ but having $a(g_1)\neq a(g_2)$ then it follows that \eqref{eqn locmax} is sharp. This is carried out in \citep*{GriJer:AEPD}.

\subsubsection{Global analysis} \label{subsubsec global}
Since the emphasis in this article is on methods common to the fields mentioned in the title, and in order to avoid having to introduce a lot of notation, we only mention one recent result from global analysis.

We need to recall the definition of the determinant of a self-adjoint, non-negative elliptic operator $P$, for example the Laplacian. Let $\lambda_1\leq\lambda_2\leq\dots$ be the positive eigenvalues of $P$. Set $\zeta_P(s)=\sum_{j=1}^\infty \lambda_j^{-s}$. This converges and is holomorphic for $s\in\C$ with sufficiently large real part, and can be shown to have a meromorphic continuation to $s\in\C$ which is regular at $s=0$. Then define $$\det P := e^{-\zeta_P'(0)}.$$ This is motivated by the fact that the analogous formula is true for invertible linear operators on a finite-dimensional space, with the usual determinant.

Consider a tube of the form \eqref{eqn Omega3} with $I=[-1,1]$ and assume it has smooth boundary. We state the result in terms of the corresponding long cylinder $X^N$ defined by \eqref{eqn def rescaling}. This consists of the ends $X_l,X_r$ and a cylinder $[-N,N]\times Y$ connecting them. \cite{MueMue:RDLTOASRD} analyze (among other things) the asymptotics of $\det P_{X^N}$ as $N\to \infty$, where $P_{X^N}$ is a Laplace type operator of product type.  Laplace type means that it is a self-adjoint differential operator whose leading part is the Laplacian, and product type means that it equals $-\frac{d^2}{d\xbar^2}+ P_Y$ over the cylindrical part $[-N,N]\times Y$, for a Laplace type operator $P_Y$ on $Y$.\footnote{More generally, $P_{X^N}$ is allowed to be a matrix-valued operator (or system of operators); that is, a vector bundle $E$ over $X^N$ is given, and $P_{X^N}$ acts on the sections of $E$. $E$ and $P_{X^N}$ are assumed to be of 'product type' over the cylindrical part of $X^N$; that is, $E$ is the pull-back of a bundle $E_Y$ over $Y$ and $P_{X^N}=-\frac{d^2}{d\xbar^2}+ P_Y$ for a Laplace type operator $P_Y$ acting on sections of $E_Y$ over $Y$. Laplace type means that, with respect to a (product type) hermitean structure on $E$, the leading part of $P_{X^N}$, when representing it in terms of an orthonormal frame of $E$, is $\Delta_{X^N}$ times the identity on $E$.
}
If $X^N$ has boundary then Dirichlet boundary conditions are imposed.

Let $X_l^\infty = \left[X_l\cup ([0,\infty)\times Y)\right]/j_l$, i.e.\ the left end with a half infinite cylinder attached at the boundary $Y$, and define $X_r^\infty$ similarly. These are the two components of the limit space introduced below in \eqref{eqn def Xinfty}. Since $P_{X^N}$ is product type over the cylinder part, it induces operators $P_{X_l^\infty}, P_{\Xinfty_r}$ on $\Xinfty_l,\Xinfty_r$.
Assuming that $P_Y$, $P_{X_l^\infty}$ and $ P_{\Xinfty_r}$ are injective (for example, $P=\Delta+1$), \citet*[Theorem 1.2]{MueMue:RDLTOASRD} prove
$$
\det P_{X^N} \sim e^{-AN}B,\quad N\to\infty
$$
for certain constants $A,B$. $A$ is calculated in terms of the function $\zeta_{P_Y}$, and $B$ in terms of similar quantities (relative determinants) on $\Xinfty_l, \Xinfty_r$, and $\det P_Y$.

There has been a lot of work in similar directions, see the references in \cite{MueMue:RDLTOASRD}. As examples we mention \cite{ParWoj:ADZDDNO}, \cite{LoyPar:DZDLMWCE}, \cite{Lee:AEZDDLCMWB}. An analysis of analytic torsion, involving similar issues, was done by \cite{Has:ASAT}.

\smallskip

{\em The basic principle and spectral invariants:}
$\det P$, like most spectral invariants occuring in global analysis, depends on {\em all} eigenvalues of $P$ simultaneously. Therefore one might expect that the basic principle, which deals with $\lambda_k(\varepsilon)$ for fixed $k$, is not applicable. However, in singular limit problems one can often separate the influence of 'small' eigenvalues from the influence of 'large' eigenvalues. In the analysis of the former the basic principle does apply, and therefore some of the methods in Section \ref{sec methods} are common to fixed $k$ and spectral invariant analysis. For the large eigenvalues one uses other techniques, for example heat kernel asymptotics as $t\to 0$.

\subsection{Graphs as limits}\label{subsec graphs}
In this section we consider tubes modeled on a finite metric graph $G=\Omega_0$. We use the notation introduced  in Section \ref{sec problem}. We first introduce quantum graphs and then  consider the following situations:
\begin{itemize}
\item
General boundary conditions, constant cross section along all edges;\footnote{That is, the situation described in  Section \ref{sec problem} with all $Y_e$ the same. The case of different $Y_e$ is also treated by \citet*{Gri:SGNS}. The main difference is that edges $e$ for which the lowest eigenvalue $\nu_0(Y_e)$ is bigger than $\min_{e'\in E}\nu_0(Y_{e'})$ should be deleted in the quantum graph of the main result in \ref{subsubsec graphs general}.}
\item
Neumann boundary condition, cross section may vary;
\end{itemize}
Thin tubes modeled on graphs have been studied extensively in the last decade. The case of Neumann boundary conditions (or empty boundary) has been well-understood for some time. For Dirichlet or more general boundary conditions, decisive progress was made only recently. We will explain some of these works. First, we introduce the limit operators on the limit graph.

\subsubsection{Quantum graphs}
What should one expect for the eigenvalues and eigenfunctions on a tube around a graph?

In the two situation mentioned above, the first cross-sectional eigenvalue is constant (say $\nu_0$) along each edge ($\nu_0=0$ in the Neumann case, with the constant function as eigenfunction; $\nu_0$ possibly positive but constant in general).
From the basic principle, see Section \ref{subsec basic principle}, one expects the $k$th eigenfunction $u_{k,\varepsilon}$ on $\Omega_\varepsilon$ to be, in the tube around any edge $e$, close to a product $u_{k}^e(x) \varphi_0(y/\varepsilon)$,
where $u_{k}^e$ is an eigenfunction of $-\frac{d^2}{dx^2}$ on the edge. Also, the eigenvalue should be approximately $\varepsilon^{-2}\nu_0 + \mu_k$, where $\mu_k$ is the eigenvalue corresponding to $u_{k}^e$. Now it must be the same $\mu_k$ for all edges $e$, since $u_{k,\varepsilon}$ is an eigenfunction on all of $\Omega_\varepsilon$.
If one knew how the $u_{k}^e$ fit together at the vertices then we could determine the $\mu_k$. Without such information the $\mu_k$ cannot be determined. This 'fitting together' is made precise in the concept of a quantum graph.

Thus, a {\em quantum graph} is, per definition, a metric graph together with a self-adjoint realization of the operator $-\frac{d^2}{dx^2}$ along each edge; such a realization is given by boundary conditions at the vertices.\footnote{One may also consider more general quantum graphs, for example allowing a magnetic or an electric field.}
One may characterize all possible such boundary conditions, see \citet*{KosSch:KRQW}. We will only encounter a special class of boundary conditions described as follows:

For each vertex $v$ denote the set of adjacent edges by $E(v)=\{e\in E:\ e\sim v\}.$
Suppose for each $v$ a linear subspace $$W_v\subset \C^{E(v)}$$
is given. This determines a boundary condition as follows: Suppose a $C^1$ function $u^e$ is given on each edge $e$ (including its endpoints). We say $u=(u^e)_{e\in E}$ satisfies the boundary condition given by the $W_v$ iff for all $v$
\begin{align}
\label{eqn bc1}
(u^e(v))_{e\sim v} &\in W_v\\
\label{eqn bc2}
(\partial_n u^e(v))_{e\sim v} &\in W_v^{\perp}.
\end{align}
Here, $u^e$ is the restriction of $u$ to the edge $e$, $\partial_n u^e(v)$ is the derivative of $u^e$ at $v$ pointing inside $e$, and $W_v^\perp$ is the orthogonal complement of $W_v$. Two special cases are especially important:

{\em Kirchhoff boundary conditions:} All $W_v=\Span\{(1,\dots,1)\}$, that is, $u$ is continuous at the vertices and the sum of derivatives at each vertex is zero.

{\em Decoupled Dirichlet conditions:} All $W_v=0$, that is, $u^e$ is zero at both ends of $e$, for each edge $e$.

For more information on quantum graphs, see the recent surveys \citep{Kuc:QGISBS} and \citep{GnuSmi:QGAQCUSS}.

\subsubsection{General boundary conditions, constant cross section}\label{subsubsec graphs general}
We assume that boundary conditions on $\partial\Omega_\varepsilon$ are arbitrary except that they respect the product decomposition $[0,l(e)]\times \varepsilon Y_e$ along each edge. That is, the boundary of each cross section $Y_e$
is divided into two parts $\calD,\calN$ and the parts of the vertex neighborhoods $X_v$ where no edges are attached are decomposed similarly. This induces an analogous decomposition of $\partial\Omega_\varepsilon$. Then we impose Dirichlet conditions on the $\calD$ part and Neumann conditions on the $\calN$ part.\footnote{Robin boundary conditions may also be considered, see \citep{Gri:SGNS}.}

\begin{figure}
\includegraphics[scale=0.8]{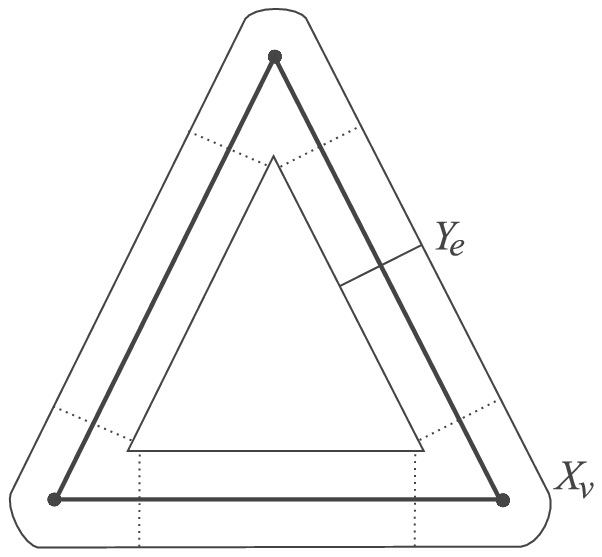} \hfill \includegraphics[scale=0.6]{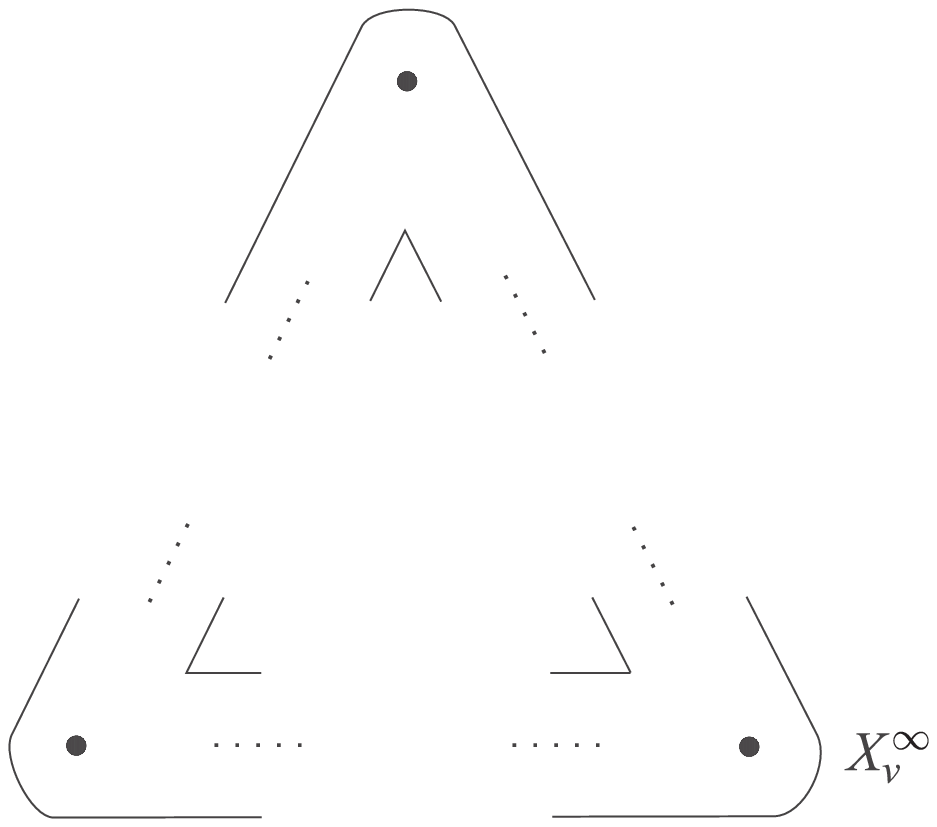}
\caption{Rescaled graph neighborhoods: $X^N$ and $\Xinfty$}
\label{figure dreieck}
\end{figure}

It turns out that the expectation stemming from the basic principle and applied in the previous section is {\em almost} correct, and that the quantum graph (i.e.\ the boundary conditions at the vertices) can be determined in terms of certain scattering data. To state this, we define the limit $X^\infty$ of the rescaled tube \eqref{eqn def rescaling} as\footnote{\label{fn Xinfty def} In the case of \eqref{eqn def phi GriJer} this would be the union of \eqref{eqn rescaled limit GriJer} and a half strip $[0,1]\times [0,\infty)$, so the notation is not consistent with \eqref{eqn rescaled limit GriJer}. However, the two ways of looking at the 'limit' of $X^N$ are equivalent, see also footnote \ref{fn scatt matrix}.}

\begin{equation}
\label{eqn def Xinfty}
X^\infty := \bigcup_v X_v^\infty \quad \text{ (disjoint union)}
\end{equation}
where the 'star' $\Xinfty_v$ of a vertex $v$ of $G$ is obtained by attaching a half infinite cylinder  $[0,\infty)\times Y_e$ to the vertex manifold $X_v$, for each edge $e$ incident to $v$, see Figure \ref{figure dreieck}.
Fix a vertex $v$.
The Laplacian $-\Delta_{X^\infty_v}$ (with  boundary conditions stemming from those on $\Omega_\varepsilon$) has continuous spectrum $[\nu_0,\infty)$, and possibly $L^2$-eigenvalues. Standard scattering theory produces for each $\lambda\in[\nu_0,\infty)$ a unitary scattering matrix
\begin{equation}
\label{eqn def scatt matrix}
S_v(\lambda) : \C^{E(v)}\to \C^{E(v)}
\end{equation}
describing the generalized eigenfunctions for the spectral value $\lambda$, i.e. the solutions $\ubar$ of
\begin{equation}
\label{eqn gen ef}
\Delta \ubar = \lambda\ubar, \quad\exists N\ |\ubar(\xbar,\ybar)| =O(\xbar^N)\text{ for }\xbar\to\infty.
\end{equation}
For the determination of the quantum graph only the scattering matrix at $\lambda=\nu_0$ is needed.
Using separation of variables it is easy to see that any solution of \eqref{eqn gen ef} with $\lambda=\nu_0$ is, along each edge cylinder $[0,\infty)\times Y_e$, asymptotically equal to either $w_e \varphi_0(\ybar)$ or $w_e \xbar \varphi_0(\ybar)$ for some $w_e\in \C$, where $\xbar$ is the coordinate along the cylinder axis. Then $S_v(\nu_0)$ is characterized by  the property that it is a unitary involution, i.e.\ it has only eigenvalues $\pm 1$ with mutually orthogonal eigenspaces, and that its $+1$-eigenspace is given by the bounded solutions, i.e.
\begin{align}
\label{eqn bdd soln}
S_v(\nu_0)w = w &\Longleftrightarrow \text{ there is a generalized eigenfunction } \ubar(\xbar,\ybar)\sim w\varphi_0(\ybar)
\end{align}
Here, if $w=(w_e)_{e\in E(v)} $ then $w\varphi_0$ is the function that equals $w_e\varphi_0$ on the edge $e$.

\smallskip
\noindent{\bf Main result on spectral convergence for tubes modeled on graphs:} {\em
Let $\mu_k$ be the eigenvalues, counted with multiplicity, of the quantum graph given by the metric graph $G$ with boundary conditions \eqref{eqn bc1}, \eqref{eqn bc2} where $W_v$ is the $+1$-eigenspace of $S_v(\nu_0)$, for each vertex $v$. Also, let $\tau_1\leq\tau_2\leq \dots\leq\tau_D$ be the $L^2$-eigenvalues $\leq\nu_0$ of $\Delta_{\Xinfty}$.
Then the eigenvalues on $\Omega_\varepsilon$ satisfy, for fixed $k$ and $\varepsilon\to 0$, and for some constant $c>0$,
\begin{align}
\label{eqnlowev Gepsilon}
\lambda_k(\varepsilon) &= \varepsilon^{-2} \tau_k + O(e^{-c/\varepsilon}),&& k=1,\dots,D\\
\label{eqnhighev Gepsilon}
\lambda_k(\varepsilon) &= \varepsilon^{-2} \nu_0 + \mu_{k-D} + O(\varepsilon),&& k>D
\end{align}
}
\smallskip

This result has features distinct from the special case of Neumann boundary conditions, which had been known before and which will be discussed below:
\begin{itemize}
\item The eigenvalues \eqref{eqnlowev Gepsilon} lie below those expected from the 'basic principle'. Their presence (i.e.\ $D>0$) is fairly typical in the case of Dirichlet boundary conditions, as is well-known in  waveguide theory, see for example \citep{SchRavWyl:QBSCUSCW}, \citep{AviBesGirMan:QBSOG} and \citep{DucExn:CIBSQWTTD}.
\item
The boundary conditions at the vertices of $G$ are determined by scattering data.
\end{itemize}
We will see in \ref{subsubsec graphs Neumann} that the behavior of the $\lambda_k(\varepsilon)$ in the Neumann problem is determined by the combinatorics of $G$ alone. These features show that this is not the case for the Dirichlet problem; here the particular geometry of the edge and vertex neighborhoods matter also.

This result is proved by \citet*{Gri:SGNS}, and partially by \citet*{MolVai:LONTFSNT}, by a different method (see Section \ref{sec methods}). \citet*{MolVai:LONTFSNT} also admit graphs with infinitely long edges (see also \citet*{MolVai:SSNTFSDA}, and \cite*{MolVai:TNTFQGESM} for a simpler case).
\citet*{Gri:SGNS} derives a full asymptotic development of $\lambda_k$ for all $k$ (cf.\ \eqref{eqn GriJer lambdak} and the comments before and after it for reasons to consider at least a third term), even uniformly for $k\leq C\varepsilon^{-1}$, for some constant $C$,\footnote{this corresponds to analyzing the full interval $[\nu_0,\nu_1)$ in the long cylinder problem} and obtains asymptotics for the eigenfunctions as well.
If the vertex neighborhoods are 'small' in the sense that their Dirichlet eigenvalues (however, with Neumann condition at the parts of the boundary where the edges are attached) are bigger than $\nu_0$ then $D=0$ and all $W_v=0$. In this case, the theorem was proved by \citet*{Pos:BQWGDBCDC}.
\smallskip

\noindent{\em A different model:}
Although the emphasis is on bounded tubes in this article, we mention interesting recent work on certain unbounded tubes modeled on the graph $\Omega_0$ having two infinitely long edges joined at one vertex, with Dirichlet boundary conditions: Given a non-intersecting smooth curve $\gamma:\R\to\R^2$ with compactly supported curvature $\kappa$ (so it consists of two straight half lines outside a compact set) define $\Omega_\varepsilon$ as the $\varepsilon$-neighborhood of the curve $t\mapsto\varepsilon^r\gamma(t)$, for some $r<1$. Thus, the curved portion (corresponding to the vertex neighborhoods in the graph model) has length $\sim \varepsilon^r$, and the same thickness $\varepsilon$ as the rest of the tube. \citet*{AlbCacFin:CSLTQW} prove, assuming $r<2/5$, convergence of the resolvent (in a certain sense) to the resolvent of the operator $-d^2/dx^2$ on $\Omega_0$, with one of two possible boundary conditions at the vertex (either Dirichlet decoupled or with a certain coupling). Which boundary condition occurs depends on whether the operator $-\frac{d^2}{dx^2}-\frac14\kappa^2$ has a bounded generalized eigenfunction at the bottom (zero) of its continuous spectrum. This is very similar to the situation in the main result above: The limit operator is the graph Laplacian, with boundary condition at the vertex determined by the rescaling of $\Omega_\varepsilon$ which makes the vertex neighborhood of a fixed 'size'.
\citet*{CacExn:NECDNSCBW} discuss a slightly more general model (with the same $\Omega_0$) in which they can produce more general boundary conditions at the vertex.

\subsubsection{Neumann boundary conditions}\label{subsubsec graphs Neumann}
In the case of Neumann boundary conditions on all of $\partial \Omega_\varepsilon$ it is easy to see that, in the main result of the previous section, $D=0$ and the quantum graph boundary conditions are Kirchhoff at each vertex. Also, $\nu_0=0$, and therefore the 'interesting' behavior occurs already in the leading term. The result $\lambda_k(\varepsilon)\to \mu_k$ $(\varepsilon\to 0)$ is then accessible by more elementary methods (quadratic forms) and is also much more stable than the result for general boundary conditions, for example under bending of the edges or varying the thickness of the tube along the edges. Such results\footnote{often with weaker error terms than $O(\varepsilon)$} were obtained  much earlier.

\citet*{Col:MPVPNNL} seems to be the first who proved this result (for straight tubes of constant width), but it  seems that the paper was not known in much of the mathematical physics community until recently.

\citet*{FreWen:DPGAP} take up the problem from a probabilistic perspective; they do not consider the spectral problem but analyze the equation $(\Delta + c)u=g$ on $\Omega_\varepsilon$ (a neighborhood in $\R^n$ of $\Omega_0$) for given continuous functions $c>0$ and $g$ and show that the solution $u$ converges for $\varepsilon\to 0$ to the solution of the corresponding equation on $\Omega_0$, with Kirchhoff boundary conditions at the vertices. This is, in essence, strong resolvent convergence.
\citet*{Sai:LENLSD} obtains a similar result for curved tubes in the plane of variable thickness.

\citet*{RubSch:VPMCTSI} prove the spectral convergence for curved tubes in the plane of variable thickness. They also allow a magnetic field. Their method was simplified and the result generalized to quite arbitrary vertex neighborhoods (of 'roughly the size' $\varepsilon$) by \citet*{KucZen:CSMSCG}, and extended to vertex neighborhoods that shrink more slowly (at a rate $\varepsilon^r$, $r<1$) by \citet*{KucZen:ASNLTD}. Depending on $r$ this may give different boundary conditions at the vertices.  \citet*{ExnPos:CSGLTM,ExnPos:QNMG} prove similar results for straight tubes of constant cross section in the manifold setting.

\section{Methods} \label{sec methods}
What makes thin tube problems difficult? There are two sources of difficulty when analyzing spectral problems of the Laplacian on $\Omega_\varepsilon$, uniformly as $\varepsilon\to 0$:
\begin{itemize}
\item
In a tube without ends, one has to deal with two scales: The transversal direction scales like $\varepsilon$ while the longitudinal direction is fixed. This is closely related to so-called semi-classical problems, i.e.\ the analysis of the operator \eqref{eqn schroedinger} as $\varepsilon\to 0$.
\item
In a tube with ends one has the additional difficulty that different parts of $\Omega_\varepsilon$ scale in different ways: the ends scale like $\varepsilon$ in {\em all} directions, the tubes around the edges only in the transversal direction.
\end{itemize}

There are various ways to deal with these problems. They differ in complexity, flexibility and precision. We mention some of them. They will be explained in more detail below.
\begin{itemize}
\item{\em Quadratic forms:} This is the simplest and most flexible tool. It is quite robust in low regularity situations. For example, it allows a fairly direct and uniform treatment of tubes around graphs where the vertex neighborhoods scale at powers $\varepsilon^r$ with $r\neq 1$, or where they are not given as scalings of a fixed model.

However, the quadratic form technique does not give optimal results in many situations, especially when one needs to understand lower order asymptotic terms.
For example, it has not been successful for the Dirichlet problem on tubes with ends (or, more generally, tubes modeled on graphs).\footnote{Another, classical example for this weakness of quadratic forms occurs with  the Weyl formula for the asymptotic distribution of eigenvalues of the Laplacian on a bounded domain $\Omega\subset\R^n$ with $C^2$ boundary: Let $N(\lambda)$ denote the number of eigenvalues $\leq\lambda$. With Dirichlet-Neumann bracketing, a quadratic forms technique, one can prove $N(\lambda)=C_n\lambda^{n/2} + R(\lambda)$ with an error estimate $R(\lambda) = O(\lambda^{(n-1)/2}\log\lambda)$, see \cite{CouHil:MMP}, for example. The sharper estimate  $R(\lambda)=O(\lambda^{(n-1)/2})$ holds but apparently cannot be proven using quadratic forms. There are various proofs of it using resolvent or heat kernel analysis.
See \cite{Gri:UBELCMWB} for further references.
On the other hand, the Dirichlet-Neumann bracketing still works for domains with very rough, for example fractal, boundary, see \citep{Lap:FDISPEOPRWBC}.
}
\item{\em Matched asymptotic expansions:} This is a classical tool for problems involving regions that scale in different ways. It can give very precise information. A technical problem is to show that the solutions (eigenfunctions) constructed by this method are {\em all} the solutions.
\item{\em Resolvent analysis:}  A good understanding of the resolvent of an operator yields information on its spectral data and also much more. Rough information on the resolvent may sometimes be obtained using quadratic forms. A precise description of the resolvent is rather complicated since in addition to the parameter $\varepsilon$ one has to analyze the dependence on the resolvent parameter $\lambda$. There are highly sophisticated techniques to do this, and when applicable they yield the most precise information. The pseudodifferential technique described below allows to treat tubes whose cross section varies mildly.
\item{\em Matching of scattering solutions:} This is a technique especially adapted to thin tube problems with constant cross section. It is rather direct and yields very precise information. As for matched asymptotics, one has to work a little harder to show that one obtains {\em all} eigenfunctions with this method.
\end{itemize}

\subsection{Quadratic forms}
The simplest method to obtain eigenvalue comparison results is to compare the quadratic forms corresponding to the operators involved: If $P$ is a  self-adjoint operator  with domain $\Dom(P)$ and quadratic form $q_P(u)=\langle u, Pu\rangle$ then the $k$th eigenvalue of $P$ is (assuming $P$ has discrete spectrum and is bounded from below), by the variational principle,
\begin{equation}
\label{eqn variation}
\lambda_k(P) = \inf_{L_k\subset\Dom(P)} \sup_{u\in L_k\setminus 0} \frac{q_P(u)}{\|u\|^2},
\end{equation}
where $L_k$ runs over $k$-dimensional subspaces.
This immediately implies that if $q_{P'}\leq q_P$ then $\lambda_k(P')\leq \lambda_k(P)$ for all $k$, and more generally:
\begin{quote}
{\bf Comparison principle:} If $J:\Dom(P)\to\Dom(P')$ is a linear map that decreases the norm at most very little and increases the quadratic form at  most very little (i.e.\ $q_{P'}(Ju)\leq q_P(u)+$ small) then $\lambda_k(P')$ can be at most slightly bigger than $\lambda_k(P)$.
\end{quote}
Here, 'very little' must be measured in an appropriate norm (e.g. the $H^1$ Sobolev norm if $P,P'$ are second order differential operators).\footnote{See, for example, \citet*{ExnPos:QNMG} for a precise statement, or \citet*{Col:MPVPNNL} for a different version.}

This principle allows for $P,P'$ to live in different Hilbert spaces and is therefore suitable for our problem, since we want to compare $\Delta_\varepsilon$, defined in $L^2(\Omega_\varepsilon)$, with $\Delta_0$, defined in $L^2(\Omega_0)$.
The quadratic forms in our context are simply $q_{\Delta_\varepsilon}(u)=\int_{\Omega_\varepsilon} |\grad u|^2$.

We explain how the comparison principle yields the eigenvalue convergence in the case of thin tubes around a graph, with Neumann boundary conditions, see \ref{subsubsec graphs Neumann} (but with a worse error $O(\sqrt\varepsilon)$ in \eqref{eqnhighev Gepsilon}; we loosely follow \citep{ExnPos:QNMG} which we refer to  for more details, and to \citep{ExnPos:CSGLTM} for the full proof along these lines). We need to find operators $J:L^2(\Omega_0)\to L^2(\Omega_\varepsilon)$ and $J':L^2(\Omega_\varepsilon)\to L^2(\Omega_0)$ respecting the domains of $\Delta$ which yield good comparison estimates of norms and quadratic forms as needed in the comparison principle. Good estimates may be expected if $J$, $J'$ reflect the expectation that eigenfunctions on $\Omega_\varepsilon$ with fixed $k$ and $\varepsilon$ small should be approximately constant in the transversal direction.

To produce from a function $f$ on the graph $\Omega_0$ a function $J(f)$ on $\Omega_\varepsilon$ is easy: Simply extend $f$ to be constant transversal to each edge and constant in each vertex neighborhood, then multiply by $\varepsilon^{-m/2}\vol(Y)^{-1/2}$ to make up for the area of the transversal slice. It is then easy to obtain the needed estimates, and this yields $\lambda_k(\Delta_\varepsilon) \leq \lambda_k(\Delta_0)+ O(\varepsilon)$.

The opposite direction is a little trickier: Given a function $u$ on $\Omega_\varepsilon$, it is natural to define the function $J'(u)$ on $\Omega_0$ along the edges by taking transversal averages, and at each vertex by taking an average of $u$ over the vertex neighborhood (and then multiplying by a suitable scaling factor). However, the resulting function may be discontinuous at the vertices, so some kind of 'fudging' has to be done. For this, various procedures have  been proposed in the literature, the simplest one being the use of a smooth transition function\footnote{For example, if $f_0$ is a given continuous function on $[0,l]$ and $a\in\R$ is close to $f_0(0)$ then one may produce a function on $[0,l]$ which is equal to $a$ at zero and close to $f_0$ everywhere by choosing a smooth $\rho:[0,l]\to[0,1]$ which equals one at $x=0$ and zero for $x\geq l/2$, and setting $f(x) = f_0(x) + \rho(x)(a-f_0(0))$.}. Then the estimate $q_0(J'(u))\leq q_\varepsilon(u)+ $'small' is fairly straightforward. The hardest part, however, is to show that $\|J'(u)\|_{L^2(\Omega_0)}$ cannot be much smaller than $\|u\|_{L^2(\Omega_\varepsilon)}$. This is done in two steps:
\begin{enumerate}
\item Consider a transversal slice to a point of an edge. Then the deviation (in $L^2$ of the slice, times the scaling factor) of $u$ from its average over the slice  is $O(\varepsilon)\|u\|_{H^1(\text{slice})}$ by the Poincar{\'e} inequality\footnote{The Poincar{\'e} inequality says that $\int_Y |u-\ave_Y u|^2 \leq C_Y \int_Y |\grad u|^2$ where $\ave_Y u$ is the average of $u$ over $Y$ and $C_Y$ is a constant depending on $Y$ (but not on $u$). A change of variable gives the scaling law $C_{\varepsilon Y} = \varepsilon^2 C_Y$.}. Integrating along the edges yields a lower bound of  $\|J'(u)\|_{L^2(\Omega_0)}$ by the $L^2$-norm of $u$ over the edge tubes (minus a small error $O(\varepsilon)\|u\|_{H^1(\Omega_\varepsilon)}$).
\item
Bound the $L^2$-norm of $u$ over the vertex neighborhoods by a small multiple of its $L^2$-norm over the edge tubes (plus small $\|u\|_{H^1}$-contributions).
\end{enumerate}
Together, this yields the bound $\|J'(u)\|_{L^2(\Omega_0)}\geq \|u\|_{L^2(\Omega_\varepsilon)}-O(\sqrt\varepsilon)\|u\|_{H^1(\Omega_\varepsilon)}$, hence
$\lambda_k(\Delta_0) \leq \lambda_k(\Delta_\varepsilon)+ O(\sqrt\varepsilon)$, which was to be shown.

It is in Step 2 that the proof hinges most essentially  on the fact that Neumann boundary conditions are considered. It can be still salvaged for Dirichlet conditions provided the vertex neighborhoods are 'small', as observed by \citet*{Pos:BQWGDBCDC}. In general however, a bound as needed here cannot hold, due to the existence of $L^2$-solutions to the associated scattering problem, see footnote \ref{fn L2 solutions}. In fact, $\Delta_\varepsilon$ will then have eigenvalues that do not come from any graph eigenvalues (those in \eqref{eqnlowev Gepsilon}).

The quadratic forms technique was also used by \citet*{Col:MPVPNNL}, \citet*{RubSch:VPMCTSI} and \citet*{KucZen:CSMSCG,KucZen:ASNLTD}.

\subsection{Resolvents, I}\label{subsec resI}
The resolvent $R_P(z)=(P-z)^{-1}$ of an operator $P$, defined for $z\in\rho(P)=\C\setminus\sigma(P)$, is a very powerful object in the study of spectral degeneration. If one can describe it sufficiently well, then one may deduce
information about the eigenvalues and also the eigenvectors: if $P$ is self-adjoint with discrete spectrum and $\gamma$ is a simple closed curve in $\rho(P)$ then
\begin{equation}
\label{eqn contour integrals}
\frac1{2\pi i}\operatorname{trace}\int_\gamma z R_P(z)\, dz,\quad \frac1{2\pi i}\int_\gamma R_P(z)\, dz
\end{equation}
are the sum of the eigenvalues of $P$ inside $\gamma$ and the orthogonal projection to the corresponding spectral subspace, respectively. Therefore, if we have a family $P_\varepsilon$  and $\gamma\subset\rho(P_\varepsilon)$ for each $\varepsilon$ then the $\varepsilon$-dependence for these objects can be deduced easily from that of $R_{P_\varepsilon}$. In the case of regular perturbations this was studied in detail by \citet*{Rel:SSI}, see also \citet*{Kat:PTLO}.

Denote by $R_\varepsilon$ the resolvent of the Laplacian on $\Omega_\varepsilon$ (in any one of the situations discussed above) and $R_0$ the resolvent of the (expected) 'limit' operator on $\Omega_0$.
Since degeneration of a tube to a one-dimensional limit object is a singular perturbation, the behavior of $R_{P_\varepsilon}$ as $\varepsilon\to 0$ is more complicated than for a regular perturbation.

 There are two approaches to understanding $R_\varepsilon$:
\begin{enumerate}
\item Compare $R_\varepsilon$ and $R_0$ directly, or
\item construct an approximation (parametrix) for $R_\varepsilon$,\footnote{this may, and probably will, involve some information coming from $R_0$} then show that it is actually close to $R_\varepsilon$.
\end{enumerate}

{\em Direct resolvent comparison:}
In this approach, an initial difficulty is that $R_\varepsilon$ and $R_0$ act on different Hilbert spaces.
In the case of a tube without ends (see \ref{subsubsec varthick I} for notation), one may write $L^2(\Omega_\varepsilon)=L_0\oplus L_{\geq 1}$
where $L_0$ corresponds to the first transversal mode and $L_{\geq 1}$ to the higher ones. $L_0$ may be identified with $L^2(I)=L^2(\Omega_0)$, and therefore $R_0$ may be regarded as operator on $L_0$, and extending it by zero on $L_{\geq 1}$ one obtains an operator on $L^2(\Omega_\varepsilon)$ which may be compared with $R_\varepsilon$. This is the approach taken by \citet*{FriSol:SDLNS,FriSol:SDLNSII}, except that they subtract the divergent term $\varepsilon^{-2}\pi^2$ from $\Delta_{\Omega_\varepsilon}$ and \eqref{eqn schroedinger} before taking inverses.\footnote{Here, only inverses are needed, not the values of the resolvent at other spectral parameters. The closeness of eigenvalues is deduced from the simple estimate $|\mu_j(A)-\mu_j(B)|\leq \|A-B\|$ for compact selfadjoint operators, applied to these inverses. This is simpler than \eqref{eqn contour integrals}. } The closeness of these operators is proved by a quadratic form comparison.

In the case of tubes with ends, or more generally graphs, a similar program was carried out by \citet*{Pos:SCQODS} (see also \citet*[section 4]{ExnPos:QNMG}) in the case of Neumann boundary conditions.
However, so far this approach has been unsuccessful when dealing with the Dirichlet problem for tubes with ends, or more generally tubes modeled on graphs.
\medskip

{\em Resolvent construction:}
\citet*{MolVai:SSNTFSDA,MolVai:LONTFSNT} use scattering solutions (see \eqref{eqn gen ef form} below) to construct a candidate approximation for the integral kernel of $R_\varepsilon$ for a tube modeled on a graph, with Dirichlet boundary conditions. They outline an argument implying that this approximates $R_\varepsilon$ to order $O(\varepsilon)$ in operator norm.
A construction of $R_\varepsilon$ to all orders was given by \citet*{HasMazMel:ASAE}. The method of the latter will be explained in \ref{subsec pseudo}.
\subsection{Matched asymptotic expansions}\label{subsec matched}
This is a general technique for problems involving different regions (sometimes called regimes) that scale in different ways. It was used by \citet*{GriJer:AEPD} in the setup introduced in \ref{subsubsec sharp loc max}. The basic steps are:
\smallskip

\noindent{\em Step 1:} We would like to find an eigenvalue $\lambda(\varepsilon)$ with eigenfunction $u_\varepsilon$ on $\Omega_\varepsilon$. One assumes that the eigenvalue has a complete asymptotic expansion as $\varepsilon\to 0$
\begin{equation}
\label{eqn ev expansion}
\lambda(\varepsilon) \sim \varepsilon^{-2}\pi^2 + \sum_{i=0}^\infty \varepsilon^i\lambda^{(i)}, \quad \lambda^{(i)}\in\R.
\end{equation}
For the eigenfunction one assumes expansions of two kinds, corresponding to the two regimes of $\Omega_\varepsilon$. The first regime is the interior of the edge: For any fixed $x\in (0,1]$ and $\ybar\in[0,1]$ one assumes an expansion refining \eqref{eqn eigenfcn approx} (with $h\equiv 1,\ h_1\equiv 0$)
\begin{equation}
\label{eqn ef expansion1}
u_\varepsilon(x,\varepsilon\ybar) \sim \sum_{i=0}^\infty \varepsilon^i u^{(i)}(x,\ybar)
\end{equation}
for functions $u^{(i)}$ on $\Omegatilde_0:=(0,1)\times[0,1]$.
The second regime is the order $\varepsilon$ neighborhood of the left end: For any fixed $(\xbar,\ybar)\in \Xinfty$ one assumes an expansion
\begin{equation}
\label{eqn ef expansion2}
u_\varepsilon(\varepsilon\xbar,\varepsilon\ybar) \sim \sum_{j=0}^\infty \varepsilon^j v^{(j)}(\xbar,\ybar)
\end{equation}
for functions $v^{(j)}$ on $X^\infty$.
\smallskip

\noindent{\em Step 2:} Plugging first \eqref{eqn ev expansion}, \eqref{eqn ef expansion1}, then \eqref{eqn ev expansion}, \eqref{eqn ef expansion2} into the equation $\Delta u_\varepsilon=\lambda(\varepsilon)u_\varepsilon$ and comparing coefficients of powers of $\varepsilon$ one obtains two infinite sequences of PDEs involving the coefficients $\lambda^{(i)},u^{(i)},v^{(j)}$. Since $u_\varepsilon$ satisfies Dirichlet boundary conditions at $\partial\Omega_\varepsilon$, the $u^{(i)}$ vanish at all sides of $\Omegatilde_0$ except possibly at $x=0$, and the $v^{(j)}$ vanish at $\partial \Xinfty$. This is not enough information to solve the equations uniquely. In addition, one needs {\em matching conditions} coming from the fact that the large $\xbar$ behavior of the $v^{(j)}$ should be related to the behavior as $x\to 0$ of the $u^{(i)}$. These conditions are obtained by setting $x=\varepsilon\xbar$ in \eqref{eqn ef expansion1}, expanding each term in Taylor series around $x=0$ and comparing coefficients of $\varepsilon$ powers with \eqref{eqn ef expansion2}.

The PDEs, with these boundary and matching conditions, can now be solved recursively. For each $\lambda^{(0)}=k^2\pi^2$, $k\in\N$, they determine $\lambda^{(i)}$ for $i\geq 1$ and $u^{(i)}$, $v^{(j)}$ for $i,j\geq 0$ uniquely, up to a normalization factor for $u_\varepsilon$.
\smallskip

\noindent{\em Step 3:} One can now construct an approximate solution by 'gluing' the formal expansions \eqref{eqn ef expansion1}, \eqref{eqn ef expansion2} in a suitable way. Cutting off the series after the $M$th term one obtains $u_{\varepsilon,k,M}$, $\lambda_{k,M}(\varepsilon)$ satisfying the boundary conditions and
\begin{equation}
\label{eqn almost solution}
\Delta u_{\varepsilon,k,M} = \lambda_{k,M}(\varepsilon)u_{\varepsilon,k,M} + O(\varepsilon^M).
\end{equation}
One needs to show that this implies that $u_{\varepsilon,k,M},\lambda_{k,M}(\varepsilon)$ are close (to order $\varepsilon^M$) to an actual eigenfunction and eigenvalue. Also, one needs to show that {\em all} eigenfunctions are obtained in this way. Both are quite straightforward in the setting of domains \eqref{eqn def phi GriJer} since by domain comparison one has good a priori estimates of the eigenvalues of $\Omega_\varepsilon$. Closeness of approximate to actual eigenfunctions follows from a spectral gap argument. This works for all $M$, so it follows that $\lambda_k(\varepsilon)$ and $u_{k,\varepsilon}$ actually do have expansions as in \eqref{eqn ev expansion}, \eqref{eqn ef expansion1}, \eqref{eqn ef expansion2}.

\subsection{Resolvents, II: Pseudodifferential operators}\label{subsec pseudo}
The classical pseudodifferential operator (\PDO) calculus provides a systematic way to study the resolvent of an elliptic operator on a compact manifold $X$ (see \citep{Shu:PDOST} for example). More precisely, if $R(\lambda,p,p')$ is the integral kernel of the resolvent, where $p,p'\in X$  and $\lambda$ is the spectral parameter, then the standard \PDO\ parametrix construction gives an explicit way to find the asymptotic expansion of $R$ as $p\to p'$ ('at the diagonal'), to all orders.

In the last two decades an intensive effort has been made (still ongoing) to extend the \PDO\ calculus to settings involving either singularities (or non-compactness) of the underlying manifold or singular perturbation problems, see for example \citep{Mel:PDOCSL,Mel:DAMWC} or \citep{Schul:PDOMWS,Schul:OAWSHMWS}. Here, the goal is to obtain the asymptotics of $R$ not only at the diagonal, but also as one or both points $p,p'$ approach the singularities (or 'infinity'), or as the perturbation parameter approaches a singular limit. R.\ Melrose and his school have developed a very geometric and general approach to this problem and applied it to many kinds of singularities and degenerations. See \citep{Gri:BBC} for an introduction to these ideas and more references.

\citet*{HasMazMel:ASAE} studied from this perspective the degeneration of a thin tube in the case of closed manifolds, with an interval as limit space. We cannot explain this here in detail but only give a rough idea. The resolvent kernel is a function $R(\varepsilon,\lambda,p,p')$ on the space $\calX=\{(\varepsilon,\lambda,p,p'):\,\varepsilon\in(0,1], \lambda\not\in\R,\ p,p'\in\Omega_\varepsilon\}$, and the goal is to obtain full information of the asymptotics at the 'boundary' and at the diagonal of this space, i.e.\ as  $\varepsilon\to 0$, $\lambda\to 0$ (in this case, $\nu_0=0$), $p\to p'$ (allowing any combination of some of these limits), uniformly over all $p,p'$. This involves a large number of different regimes, similar to but more complex than the situation in \ref{subsec matched} where we considered different regimes for the eigenfunction as function of $\varepsilon$ and $p$.
The goal is accomplished in essence by a combination of the classical \PDO\ technique with the method of matched asymptotic expansions. A difficulty lies in taming the high combinatorial complexity of all the regimes and resulting matchings. This is made feasible by encoding these geometrically in a manifold with corners $\calXbar$ which is a compactification of $\calX$, whose minimal boundary faces  correspond to the regimes and boundary hypersurfaces to matchings, cf.\ \citep[section 2.5.3]{Gri:BBC}. See \citep[Figure 2]{HasMazMel:ASAE} for a nice picture of $\calXbar$. More precisely, the work proceeds by first defining a suitable space $\calXbar$ (motivated, among other things, by model calculations), then defining a class of operators (a '$\calXbar$-pseudodifferential calculus') as those whose integral kernels are distributions on $\calXbar$ having certain 'nice' (polyhomogeneous conormal) behavior at the boundary of $\calXbar$ and conormal (\PDO) singularities at the diagonal, then showing that this class is closed under composition and has various other 'good' properties, then constructing a first order parametrix for $\Delta-\lambda$ (this amounts to solving certain limit or model problems, one of which being the inversion of the Laplacian on an interval, another the scattering problem on $X^\infty$) and then improving this to any order by the standard iteration in parametrix constructions.
This can then be used to obtain full asymptotics of eigenvalues and eigenfunctions by \eqref{eqn contour integrals}.

\subsection{Matching of scattering solutions}
This very effective method was used by \cite{Gri:SGNS} for the proof of the main result in \ref{subsubsec graphs general}. We explain it in the simple case of a single edge graph with one end attached, as in \eqref{eqn def phi GriJer}.\footnote{but the special structure of the end in \eqref{eqn def phi GriJer} is inessential here, in particular $\Omega_\varepsilon$ need not be contained in the strip $0<y<\varepsilon$}

We work on the rescaled domain $X^N=\{(\xbar,\ybar):\ 0<\ybar<1,\ -g(\ybar)<\xbar<N\}$, and its limit $\Xinfty$ defined in \eqref{eqn rescaled limit GriJer}. We first need some facts from scattering theory on $\Xinfty$, extending \eqref{eqn U asymptotics}. For each $\alpha\in(0,\sqrt3 \pi)$ there is a unique bounded solution $E_\alpha$ of
\begin{equation}
\label{eqn gen ef eqn}
\Delta E_\alpha = \lambda E_\alpha\text{ on }X^\infty, \quad \lambda=\pi^2+\alpha^2,\quad E_{\alpha|\partial\Xinfty}=0,
\end{equation}
which has the form
\begin{equation}
\label{eqn gen ef form}
E_\alpha(\xbar,\ybar) = \left( e^{-i\alpha\xbar} + S(\alpha)e^{i\alpha\xbar}\right) \sin \pi\ybar + R_\alpha(\xbar,\ybar)
\end{equation}
for $\xbar\geq 0$, where $S(\alpha)\in \C$ and $R_\alpha(\xbar,\cdot)$ is in the span of the higher transversal modes $\sin\pi m\ybar$, $m\geq 2$, and consequently is exponentially decaying as $\xbar\to\infty$.\footnote{This decay follows from $R_\alpha(\xbar,\ybar)=\sum_{m=2}^\infty r_m(\xbar)\sin\pi m\ybar$ with $r_m(\xbar)=2\int_0^1 R_\alpha(\xbar,\ybar)\sin\pi m\ybar\,d\ybar$, hence bounded, and satisfying $r_m'' - \sqrt{m^2\pi^2-\lambda}r_m=0$ from \eqref{eqn gen ef eqn}, which has only exponential solutions for $\lambda<4\pi^2$.}
$S(\alpha)$ is called the scattering matrix.\footnote{\label{fn scatt matrix} In this simple example it is a number, i.e. a $1\times 1$ matrix. Note that in the general setup of \citep{Gri:SGNS} this example would be treated as the union of two half strips, one emanating from each end (vertex), see \eqref{eqn def Xinfty} and footnote \ref{fn Xinfty def}, and the scattering matrix would be the  $2\times 2$ matrix $\diag (S(\alpha),-1)$; it is easy to see from this that the two treatments are equivalent.}
It is unitary and extends to a holomorphic function of $\alpha$ in a neighborhood of $0$ in $\C$, satisfying $S(0)^2=1$. In our example $S(0)=-1$ (this means that the $+1$ eigenspace in \eqref{eqn bdd soln} is zero). Write
\begin{equation}
\label{eqn def rho}
S(\alpha)=e^{i\rho(\alpha)},\quad \rho\text{ holomorphic in }\alpha, \ \rho(0)=\pi.
\end{equation}
$\rho(\alpha)$ is real for real $\alpha$. It is called the scattering phase.\footnote{\label{fn scatt phase}The relation of \eqref{eqn gen ef form} with \eqref{eqn U asymptotics}, which is a scattering solution for $\alpha=0$, is this: Write $\rho(\alpha)=\pi+\alpha\sigma(\alpha)$ with $\sigma(0)=\rho'(0)$. Then $E_\alpha(\xbar,\ybar)=-2i e^{i\alpha\sigma(\alpha)/2} \sin \left[(\xbar+\frac{\sigma(\alpha)}2)\alpha\right] \sin\pi\ybar  + R_\alpha(\xbar,\ybar)$. Therefore the first transversal mode of $\alpha^{-1}E_\alpha$ converges, as $\alpha\to 0$, to $-2i (\xbar+\frac{\sigma(0)}2)\sin\pi\ybar$. It is easy to see that $\alpha^{-1} R_\alpha$ also converges. Therefore $U=\frac i2\lim_{\alpha\to 0} \alpha^{-1} E_\alpha$ and $a(g)=\rho'(0)/2$. }

Fix $N>0$. Suppose $\alpha$ is such that
\begin{equation}
\label{eqn alpha cond}
e^{-i\alpha N} + S(\alpha) e^{i\alpha N}= 0.
\end{equation}
Then $E_\alpha(N,\ybar)=R_\alpha(N,\ybar)=O(e^{-cN})$. Therefore, $E_\alpha$ satisfies Dirichlet boundary conditions at $\partial X^N$, up to an exponentially small error. Also, $\Delta E_\alpha=\lambda E_\alpha$ in $X^N$. From this it is easy to conclude that $\lambda$ lies within $O(e^{-cN})$ of a Dirichlet eigenvalue of $\Delta$ on ${X^N}$.\footnote{\label{fn L2 solutions}Similarly, $L^2$-eigenfunctions on $\Xinfty$ yield eigenvalues \eqref{eqnlowev Gepsilon} since they decay exponentially as $\xbar\to\infty$.
}

Therefore, we need to analyze the solutions $\alpha$ of \eqref{eqn alpha cond}. Use \eqref{eqn def rho} and introduce the variables $z=\alpha N$, $\varepsilon=N^{-1}$, then this is equivalent to
\begin{equation}
\label{eqn rho z}
\rho(\varepsilon z) + 2z = 2\pi k + \pi\quad\text{ for some }k\in\Z.
\end{equation}
By the implicit function theorem this can be solved locally for $z$, that is, for each $k$ there is a unique analytic function $\varepsilon\mapsto z_k(\varepsilon)$, defined near zero, so that $z=z_k(\varepsilon)$ satisfies \eqref{eqn rho z} and $z_k(0)=k\pi$. Then $\alpha_k(\varepsilon)=\varepsilon z_k(\varepsilon)$ solves \eqref{eqn alpha cond} with $N=\varepsilon^{-1}$, so one gets an eigenvalue $\pi^2 + \varepsilon^2 z_k(\varepsilon)^2 + O(e^{-c/\varepsilon})$ on $X^N$ and hence an eigenvalue
\begin{equation}
\label{eqn ev candidate}
\varepsilon^{-2}\pi^2 + k^2\pi^2 + \sum_{i=1}^\infty \lambda_k^{(i)} \varepsilon^i + O(e^{-c/\varepsilon})
\end{equation}
for $\Omega_\varepsilon$.\footnote{This is stronger than \eqref{eqn ev expansion} since now we know that the $\varepsilon$-series actually converges for small $\varepsilon$.} Also, $E_{\alpha_k(\varepsilon)}$ is an approximation, with error $O(e^{-c/\varepsilon})$, for the eigenfunction.

As in \ref{subsec matched} one needs to show that one has obtained {\em all} eigenvalues in this way. For general graphs, this requires a rather involved argument, which essentially rests on a priori estimates for the eigenfunctions, see  \citep{Gri:SGNS}. There are various pitfalls complicating the matter, for example the issue of multiplicities and almost multiplicities (eigenvalues can be very close together, so the spectral gap argument alluded to in \ref{subsec matched} does not work).

The method of matching of scattering solutions was also used in the global analysis context by \cite{CapLeeMil:SAEOMDILES}, \cite{Mue:EIMWB}, \cite{ParWoj:ADZDST}.


\begin{thebibliography}{71}
\providecommand{\natexlab}[1]{#1}
\providecommand{\url}[1]{\texttt{#1}}
\expandafter\ifx\csname urlstyle\endcsname\relax
  \providecommand{\doi}[1]{doi: #1}\else
  \providecommand{\doi}{doi: \begingroup \urlstyle{rm}\Url}\fi

\bibitem[Albeverio et~al.(2007)Albeverio, Cacciapuoti, and
  Finco]{AlbCacFin:CSLTQW}
S.~Albeverio, C.~Cacciapuoti, and D.~Finco.
\newblock Coupling in the singular limit of thin quantum waveguides.
\newblock \emph{Journal of Mathematical Physics}, 48\penalty0 (3):\penalty0
  032103, 2007.
\newblock \doi{10.1063/1.2710197}.
\newblock URL \url{http://link.aip.org/link/?JMP/48/032103/1}.

\bibitem[Avishai et~al.(1991)Avishai, Bessis, Giraud, and
  Mantica]{AviBesGirMan:QBSOG}
Y.~Avishai, D.~Bessis, B.~G. Giraud, and G.~Mantica.
\newblock Quantum bound states in open geometries.
\newblock \emph{Phys. Rev. B}, 44\penalty0 (15):\penalty0 8028--8034, Oct 1991.
\newblock \doi{10.1103/PhysRevB.44.8028}.

\bibitem[Berkolaiko et~al.(2006)Berkolaiko, Carlson, Fulling, and
  Kuchment]{BerCarFul:QGTA}
G.~Berkolaiko, R.~Carlson, S.~A. Fulling, and P.~Kuchment, editors.
\newblock \emph{{Quantum graphs and their applications. Proceedings of an
  AMS-IMS-SIAM joint summer research conference on quantum graphs and their
  applications, Snowbird, UT, USA, June 19--23, 2005.}}, 2006. {Contemporary
  Mathematics 415. Providence, RI: American Mathematical Society (AMS)}.

\bibitem[Bouchitt{\'e} et~al.(2007)Bouchitt{\'e}, Mascarenhas, and
  Trabucho]{BouMasTra:CTEODW}
G.~Bouchitt{\'e}, M.~Mascarenhas, and L.~Trabucho.
\newblock {On the curvature and torsion effects in one dimensional waveguides}.
\newblock \emph{ESAIM, Control Optim. Calc. Var.}, 13\penalty0 (4):\penalty0
  793--808, 2007.

\bibitem[Cacciapuoti and Exner(2007)]{CacExn:NECDNSCBW}
C.~Cacciapuoti and P.~Exner.
\newblock Nontrivial edge coupling from a {D}irichlet network squeezing: the
  case of a bent waveguide.
\newblock Preprint, arXiv:0704.2912, 2007.

\bibitem[Cappell et~al.(1996)Cappell, Lee, and Miller]{CapLeeMil:SAEOMDILES}
S.~E. Cappell, R.~Lee, and E.~Y. Miller.
\newblock {Self-adjoint elliptic operators and manifold decompositions. I: Low
  eigenmodes and stretching}.
\newblock \emph{Commun. Pure Appl. Math.}, 49\penalty0 (8):\penalty0 825--866,
  1996.

\bibitem[Cheng(1976)]{Che:ENS}
S.-Y. Cheng.
\newblock {Eigenfunctions and nodal sets}.
\newblock \emph{Comment. Math. Helv.}, 51:\penalty0 43--55, 1976.

\bibitem[Colin~de Verdi{\`e}re(1986)]{Col:MPVPNNL}
Y.~Colin~de Verdi{\`e}re.
\newblock {Sur la multiplicit{\'e} de la premi{\`e}re valeur propre non nulle
  du {L}aplacien}.
\newblock \emph{Comment. Math. Helv.}, 61:\penalty0 254--270, 1986.

\bibitem[Colin~de {V}erdi{\`e}re(1987)]{Col:CLDPFSD}
Y.~Colin~de {V}erdi{\`e}re.
\newblock {Construction de laplaciens dont une partie finie du spectre est
  donn{\'e}e}.
\newblock \emph{Ann. Sci. \'Ec. Norm. Sup\'er.}, 20:\penalty0 599--615, 1987.

\bibitem[Courant and Hilbert(1924)]{CouHil:MMP}
R.~Courant and D.~Hilbert.
\newblock \emph{Methoden der mathematischen {Physik}, Band I, II}.
\newblock Springer Verlag, Berlin, 1924.

\bibitem[Duclos and Exner(1995)]{DucExn:CIBSQWTTD}
P.~Duclos and P.~Exner.
\newblock {Curvature-induced bound states in quantum waveguides in two and
  three dimensions}.
\newblock \emph{Rev. Math. Phys.}, 7\penalty0 (1):\penalty0 73--102, 1995.

\bibitem[Ekholm et~al.(2005)Ekholm, Kovarik, and
  Krej\v{c}i\v{r}\'{i}k]{EkhKovKre:HITW}
T.~Ekholm, H.~Kovarik, and D.~Krej\v{c}i\v{r}\'{i}k.
\newblock A hardy inequality in twisted waveguides.
\newblock Preprint, arXiv:math-ph/0512050, to appear in Arch.Rat.Mech.Anal.,
  2005.

\bibitem[Exner and Post(2005)]{ExnPos:CSGLTM}
P.~Exner and O.~Post.
\newblock {Convergence of spectra of graph-like thin manifolds}.
\newblock \emph{J. Geom. Phys.}, 54\penalty0 (1):\penalty0 77--115, 2005.

\bibitem[Exner and Post(2007)]{ExnPos:QNMG}
P.~Exner and O.~Post.
\newblock Quantum networks modelled by graphs.
\newblock Preprint, arXiv:0706.0481v1, 2007.

\bibitem[Exner and {\v{S}}eba(1989{\natexlab{a}})]{ExnSeb:BSCQW}
P.~Exner and P.~{\v{S}}eba.
\newblock {Bound states in curved quantum waveguides}.
\newblock \emph{J. Math. Phys.}, 30\penalty0 (11):\penalty0 2574--2580,
  1989{\natexlab{a}}.

\bibitem[Exner and {\v{S}}eba(1989{\natexlab{b}})]{ExnSeb:ESMCOT}
P.~Exner and P.~{\v{S}}eba.
\newblock {Electrons in semiconductor microstructures: A challenge to operator
  theorists}.
\newblock {Schr\"odinger operators, standard and non-standard, Proc. Conf.,
  Dubna/USSR 1988, 79-100 (1989).}, 1989{\natexlab{b}}.

\bibitem[Freidlin and Wentzell(1993)]{FreWen:DPGAP}
M.~I. Freidlin and A.~D. Wentzell.
\newblock {Diffusion processes on graphs and the averaging principle}.
\newblock \emph{Ann. Probab.}, 21\penalty0 (4):\penalty0 2215--2245, 1993.

\bibitem[Freitas and Krej\v{c}i\v{r}\'{i}k(2007)]{FreKre:LNSTCT}
P.~Freitas and D.~Krej\v{c}i\v{r}\'{i}k.
\newblock Location of the nodal set for thin curved tubes.
\newblock To appear in Indiana Univ. Math. J., 2007.

\bibitem[Friedlander and Solomyak(2007{\natexlab{a}})]{FriSol:SDLNS}
L.~Friedlander and M.~Solomyak.
\newblock On the spectrum of the {D}irichlet {L}aplacian in a narrow strip.
\newblock Preprint, arXiv:0705.4058, 2007{\natexlab{a}}.

\bibitem[Friedlander and Solomyak(2007{\natexlab{b}})]{FriSol:SDLNSII}
L.~Friedlander and M.~Solomyak.
\newblock On the spectrum of the {D}irichlet {L}aplacian in a narrow strip, ii.
\newblock Preprint, arXiv:0710.1886, 2007{\natexlab{b}}.

\bibitem[Froese and Herbst(2001)]{FroHer:RHCCQM}
R.~Froese and I.~Herbst.
\newblock {Realizing holonomic constraints in classical and quantum mechanics}.
\newblock \emph{Commun. Math. Phys.}, 220\penalty0 (3):\penalty0 489--535,
  2001.

\bibitem[Geer and Keller(1983)]{GeeKel:ESCWST}
J.~F. Geer and J.~B. Keller.
\newblock {Eigenvalues of slender cavities and waves in slender tubes}.
\newblock \emph{J. Acoust. Soc. Am.}, 74:\penalty0 1895--1904, 1983.

\bibitem[Gnutzmann and Smilansky(2006)]{GnuSmi:QGAQCUSS}
S.~Gnutzmann and U.~Smilansky.
\newblock Quantum graphs: {A}pplications to quantum chaos and universal
  spectral statistics.
\newblock \emph{Advances in Physics}, 55\penalty0 (5-6):\penalty0 527--625,
  2006.

\bibitem[Grieser(2001)]{Gri:BBC}
D.~Grieser.
\newblock Basics of the {$b$}-calculus.
\newblock In J.~Gil, D.~Grieser, and M.~Lesch, editors, \emph{Approaches to
  Singular Analysis}, Advances in Partial Differential Equations, Basel, 2001.
  Birkh\"auser.

\bibitem[Grieser(2002)]{Gri:UBELCMWB}
D.~Grieser.
\newblock Uniform bounds for eigenfunctions of the {Laplacian} on compact
  manifolds with boundary.
\newblock \emph{Comm. PDE}, 27:\penalty0 1283--1299, 2002.

\bibitem[Grieser(2007)]{Gri:SGNS}
D.~Grieser.
\newblock Spectra of graph neighborhoods and scattering.
\newblock Preprint arXiv:0710.3405, 2007.

\bibitem[Grieser and Jerison(1996)]{GriJer:AFNLCD}
D.~Grieser and D.~Jerison.
\newblock Asymptotics of the first nodal line of a convex domain.
\newblock \emph{Invent. Math.}, 125\penalty0 (2):\penalty0 197--219, 1996.

\bibitem[Grieser and Jerison(1998)]{GriJer:SFECPD}
D.~Grieser and D.~Jerison.
\newblock The size of the first eigenfunction of a convex planar domain.
\newblock \emph{J. Am. Math. Soc.}, 11\penalty0 (1):\penalty0 41--72, 1998.

\bibitem[Grieser and Jerison(2007)]{GriJer:AEPD}
D.~Grieser and D.~Jerison.
\newblock Asymptotics of eigenfunctions on plane domains.
\newblock Preprint, arXiv:0710.3665, 2007.

\bibitem[Hassell(1998)]{Has:ASAT}
A.~Hassell.
\newblock {Analytic surgery and analytic torsion}.
\newblock \emph{Commun. Anal. Geom.}, 6\penalty0 (2):\penalty0 255--289, 1998.

\bibitem[Hassell et~al.(1995)Hassell, Mazzeo, and Melrose]{HasMazMel:ASAE}
A.~Hassell, R.~Mazzeo, and R.~B. Melrose.
\newblock {Analytic surgery and the accumulation of eigenvalues}.
\newblock \emph{Commun. Anal. Geom.}, 3\penalty0 (1):\penalty0 115--222, 1995.

\bibitem[Hassell and Zelditch(2004)]{HasZel:QEBVE}
A.~Hassell and S.~Zelditch.
\newblock {Quantum ergodicity of boundary values of eigenfunctions}.
\newblock \emph{Commun. Math. Phys.}, 248\penalty0 (1):\penalty0 119--168,
  2004.

\bibitem[Hoffmann-Ostenhof et~al.(1997)Hoffmann-Ostenhof, Hoffmann-Ostenhof,
  and Nadirashvili]{HofHofNad:NLSELRCC}
M.~Hoffmann-Ostenhof, T.~Hoffmann-Ostenhof, and N.~Nadirashvili.
\newblock The nodal line of the second eigenfunction of the {Laplacian} in
  {$R^2$} can be closed.
\newblock \emph{Duke Math. J.}, 90:\penalty0 631--640, 1997.

\bibitem[Jerison(1991)]{Jer:FNLCPD}
D.~Jerison.
\newblock {The first nodal line of a convex planar domain}.
\newblock \emph{Int. Math. Res. Not.}, 1991\penalty0 (1):\penalty0 1--5, 1991.

\bibitem[Jerison(1995{\natexlab{a}})]{Jer:DFNLCD}
D.~Jerison.
\newblock The diameter of the first nodal line of a convex domain.
\newblock \emph{Ann. Math., II. Ser.}, 141\penalty0 (1):\penalty0 1--33,
  1995{\natexlab{a}}.

\bibitem[Jerison(1995{\natexlab{b}})]{Jer:FNSCD}
D.~Jerison.
\newblock The first nodal set of a convex domain.
\newblock In C.~Fefferman et~al., editors, \emph{Essays on Fourier analysis in
  honor of Elias M. Stein. Proceedings of the Princeton conference on harmonic
  analysis held at Princeton Univ. in honor of Elias M. Stein's 60th birthday},
  volume~42 of \emph{Princeton Math. Ser.}, pages 225--249. Princeton Univ.
  Press, 1995{\natexlab{b}}.

\bibitem[Jerison(2000)]{Jer:LFNLNP}
D.~Jerison.
\newblock Locating the first nodal line in the {N}eumann problem.
\newblock \emph{Trans. Am. Math. Soc.}, 352\penalty0 (5):\penalty0 2301--2317,
  2000.

\bibitem[Kato(1976)]{Kat:PTLO}
T.~Kato.
\newblock \emph{{Perturbation theory for linear operators. 2nd ed.}}
\newblock {Grundlehren der mathematischen Wissenschaften. 132.
  Berlin-Heidelberg-New York: Springer-Verlag. }, 1976.

\bibitem[Kostrykin and Schrader(1999)]{KosSch:KRQW}
V.~Kostrykin and R.~Schrader.
\newblock {Kirchhoff's rule for quantum wires}.
\newblock \emph{J. Phys. A, Math. Gen.}, 32\penalty0 (4):\penalty0 595--630,
  1999.

\bibitem[Krej\v{c}i\v{r}\'{i}k(2008)]{Kre:TVBQW}
D.~Krej\v{c}i\v{r}\'{i}k.
\newblock Twisting versus bending in quantum waveguides.
\newblock In \emph{Analysis on Graphs and its Applications, Proceedings of
  Symposia in Pure Mathematics, P. Exner, J. Keating, P. Kuchment, T. Sunada,
  A.Teplyaev (eds.)}. AMS, 2008.

\bibitem[Kuchment(2002)]{Kuc:GMWTS}
P.~Kuchment.
\newblock {Graph models for waves in thin structures}.
\newblock \emph{Waves in Random Media}, 12:\penalty0 R1--R24, 2002.

\bibitem[Kuchment(2004)]{Kuc:QGISBS}
P.~Kuchment.
\newblock {Quantum graphs. I: Some basic structures}.
\newblock \emph{Waves Random Media}, 4:\penalty0 S107--S128, 2004.

\bibitem[Kuchment and Zeng(2001)]{KucZen:CSMSCG}
P.~Kuchment and H.~Zeng.
\newblock {Convergence of spectra of mesoscopic systems collapsing onto a
  graph}.
\newblock \emph{J. Math. Anal. Appl.}, 258\penalty0 (2):\penalty0 671--700,
  2001.

\bibitem[Kuchment and Zeng(2003)]{KucZen:ASNLTD}
P.~Kuchment and H.~Zeng.
\newblock {Asymptotics of spectra of {N}eumann {L}aplacians in thin domains}.
\newblock {Karpeshina, Yulia (ed.) et al., Advances in differential equations
  and mathematical physics. Proceedings of the 9th UAB international
  conference, University of Alabama, Birmingham, AL, USA, March 26--30, 2002.
  Providence, RI: American Mathematical Society (AMS). Contemp. Math. 327,
  199-213}, 2003.

\bibitem[Lapidus(1991)]{Lap:FDISPEOPRWBC}
M.~L. Lapidus.
\newblock {Fractal drum, inverse spectral problems for elliptic operators and a
  partial resolution of the Weyl-Berry conjecture}.
\newblock \emph{Trans. Am. Math. Soc.}, 325\penalty0 (2):\penalty0 465--529,
  1991.

\bibitem[Lee(2006)]{Lee:AEZDDLCMWB}
Y.~Lee.
\newblock {Asymptotic expansions of the zeta-determinants of Dirac Laplacians
  on a compact manifold with boundary.}
\newblock \emph{J. Geom. Anal.}, 16\penalty0 (4):\penalty0 633--660, 2006.

\bibitem[Loya and Park(2004)]{LoyPar:DZDLMWCE}
P.~Loya and J.~Park.
\newblock {Decomposition of the $\zeta$-determinant for the Laplacian on
  manifolds with cylindrical end}.
\newblock \emph{Ill. J. Math.}, 48\penalty0 (4):\penalty0 1279--1303, 2004.

\bibitem[Maz'ja et~al.(1991)Maz'ja, Nazarov, and
  Plamenevskij]{MazNazPla:ATERSGG}
V.~Maz'ja, S.~Nazarov, and P.~Plamenevskij.
\newblock \emph{{Asymptotische Theorie elliptischer Randwertaufgaben in
  singul\"ar gest\"orten Gebieten, Band I,II}}.
\newblock Akademie Verlag, Berlin, 1991.

\bibitem[Melas(1992)]{Mel:NLSELR}
A.~Melas.
\newblock On the nodal line of the second eigenfunction of the {Laplacian} in
  {$R^2$}.
\newblock \emph{J. Diff. Geom.}, 35:\penalty0 255--263, 1992.

\bibitem[Melrose(1991)]{Mel:PDOCSL}
R.~B. Melrose.
\newblock {Pseudodifferential operators, corners and singular limits}.
\newblock In \emph{Proc. Int. Congr. Math., Kyoto/Japan 1990, Vol. I}, pages
  217--234, 1991.

\bibitem[Melrose(1996)]{Mel:DAMWC}
R.~B. Melrose.
\newblock Differential analysis on manifolds with corners.
\newblock Book in preparation. http://www-math.mit.edu/$\sim$rbm/book.html,
  1996.

\bibitem[Molchanov and Vainberg(2006{\natexlab{a}})]{MolVai:SSNTFSDA}
S.~Molchanov and B.~Vainberg.
\newblock Scattering solutions in networks of thin fibers: {S}mall diameter
  asymptotics.
\newblock Preprint, arXiv:math-ph/060902, 2006{\natexlab{a}}.

\bibitem[Molchanov and Vainberg(2006{\natexlab{b}})]{MolVai:TNTFQGESM}
S.~Molchanov and B.~Vainberg.
\newblock {Transition from a network of thin fibers to the quantum graphs: an
  explicitly solvable model}.
\newblock {Berkolaiko, Gregory et al. (eds.), Quantum graphs and their
  applications. Providence, RI: American Mathematical Society (AMS).
  Contemporary Mathematics 415, 227-239 (2006).}, 2006{\natexlab{b}}.

\bibitem[Molchanov and Vainberg(2007)]{MolVai:LONTFSNT}
S.~Molchanov and B.~Vainberg.
\newblock Laplace operator in networks of thin fibers: {S}pectrum near the
  threshold.
\newblock Preprint, arXiv:0704.2795, 2007.

\bibitem[M{\"u}ller and M{\"u}ller(2006)]{MueMue:RDLTOASRD}
J.~M{\"u}ller and W.~M{\"u}ller.
\newblock {Regularized determinants of Laplace-type operators, analytic
  surgery, and relative determinants.}
\newblock \emph{Duke Math. J.}, 133\penalty0 (2):\penalty0 259--312, 2006.

\bibitem[M{\"u}ller(1994)]{Mue:EIMWB}
W.~M{\"u}ller.
\newblock {Eta invariants and manifolds with boundary}.
\newblock \emph{J. Differ. Geom.}, 40\penalty0 (2):\penalty0 311--377, 1994.

\bibitem[Park and Wojciechowski(2005)]{ParWoj:ADZDDNO}
J.~Park and K.~P. Wojciechowski.
\newblock {Adiabatic decomposition of the $\zeta$-determinant and Dirichlet to
  Neumann operator}.
\newblock \emph{J. Geom. Phys.}, 55\penalty0 (3):\penalty0 241--266, 2005.

\bibitem[Park and Wojciechowski(2006)]{ParWoj:ADZDST}
J.~Park and K.~P. Wojciechowski.
\newblock {Adiabatic decomposition of the $\zeta$-determinant and scattering
  theory}.
\newblock \emph{Mich. Math. J.}, 54\penalty0 (1):\penalty0 207--238, 2006.

\bibitem[Pavlov(2007)]{Pav:SGMOE}
B.~Pavlov.
\newblock {A star-graph model via operator extension}.
\newblock \emph{Math. Proc. Camb. Philos. Soc.}, 142\penalty0 (2):\penalty0
  365--384, 2007.

\bibitem[Payne(1973)]{Pay:TCFMEP}
L.~E. Payne.
\newblock On two conjectures in the fixed membrane eigenvalue problem.
\newblock \emph{Z. angew. Math. Phys.}, 24:\penalty0 721--729, 1973.

\bibitem[Post(2005)]{Pos:BQWGDBCDC}
O.~Post.
\newblock {Branched quantum wave guides with Dirichlet boundary conditions: the
  decoupling case}.
\newblock \emph{J. Phys. A, Math. Gen.}, 38\penalty0 (22):\penalty0 4917--4931,
  2005.

\bibitem[Post(2006)]{Pos:SCQODS}
O.~Post.
\newblock Spectral convergence of quasi-one-dimensional spaces.
\newblock \emph{Ann. Henri Poincar{\'e}}, 7:\penalty0 933--973, 2006.

\bibitem[Rellich(1936)]{Rel:SSI}
F.~Rellich.
\newblock {St{\"o}rungstheorie der {S}pektralzerlegung. I. {A}nalytische
  {S}t{\"o}rung der isolierten {P}unkteigenwerte eines beschr{\"a}nkten
  {O}perators.}
\newblock \emph{Math. Ann.}, 113:\penalty0 600--619, 1936.

\bibitem[Rubinstein and Schatzman(2001)]{RubSch:VPMCTSI}
J.~Rubinstein and M.~Schatzman.
\newblock {Variational problems on multiply connected thin strips. I: Basic
  estimates and convergence of the Laplacian spectrum}.
\newblock \emph{Arch. Ration. Mech. Anal.}, 160\penalty0 (4):\penalty0
  271--308, 2001.

\bibitem[Ruedenberg and Scherr(1953)]{RueSch:FENMCSI}
K.~Ruedenberg and C.~W. Scherr.
\newblock Free-electron network model for conjugated systems. {I}. {T}heory.
\newblock \emph{J. Chem. Phys.}, 21\penalty0 (9):\penalty0 1565--1581, 1953.

\bibitem[Saito(2000)]{Sai:LENLSD}
Y.~Saito.
\newblock {The limiting equation for {N}eumann {L}aplacians on shrinking
  domains}.
\newblock \emph{Electron. J. Differ. Equ.}, 2000.

\bibitem[Schult et~al.(1989)Schult, Ravenhall, and Wyld]{SchRavWyl:QBSCUSCW}
R.~L. Schult, D.~G. Ravenhall, and H.~W. Wyld.
\newblock Quantum bound states in a classically unbound system of crossed
  wires.
\newblock \emph{Phys. Rev. B}, 39\penalty0 (8):\penalty0 5476--5479, Mar 1989.
\newblock \doi{10.1103/PhysRevB.39.5476}.

\bibitem[Schulze(1991)]{Schul:PDOMWS}
B.-W. Schulze.
\newblock \emph{Pseudo-Differential Operators on Manifolds with Singularities},
  volume~24 of \emph{Stud. Math. Appl.}
\newblock North-Holland, Amsterdam, 1991.

\bibitem[Schulze(2001)]{Schul:OAWSHMWS}
B.-W. Schulze.
\newblock Operator algebras with symbol hierarchies on manifolds with
  singularities.
\newblock In J.~Gil, D.~Grieser, and L.~M., editors, \emph{Advances in Partial
  Differential Equations (Approaches to Singular Analysis)}, Oper. Theory Adv.
  Appl. Birkh\"auser Verlag, Basel, 2001.

\bibitem[Shubin(1987)]{Shu:PDOST}
M.~A. Shubin.
\newblock \emph{Pseudodifferential operators and spectral theory}.
\newblock Springer-Verlag, 1987.

\bibitem[Smilansky and Solomyak(2006)]{SmiSol:QGLNPW}
U.~Smilansky and M.~Solomyak.
\newblock {The quantum graph as a limit of a network of physical wires}.
\newblock {Berkolaiko, Gregory et al. (eds.), Quantum graphs and their
  applications. Providence, RI: American Mathematical Society (AMS).
  Contemporary Mathematics 415, 283-291 (2006).}, 2006.

\end{thebibliography}
\def\cprime{$'$}  \newcommand{\nop}[1]{} \newcommand{\single}[1]{#1}
  \newcommand{\SwapArgs}[2]{#2#1}
  \newcommand{\translationof}{\iflanguage{german}{Englische "Ubersetzung von }
  {English Translation of }}
  \newcommand{\submitted}{\iflanguage{german}{Eingereicht }{Submitted }}
  \newcommand{\submittedto}{\iflanguage{german}{Eingereicht bei }{Submitted to
  }} \newcommand{\privcomm}{\iflanguage{german}{Pers\"onliche
  Mitteilung}{Private communication}}
\addtocontents{toc}{\SkipTocEntry}  

\end{document}